\documentclass[11pt,a4paper]{article}

\usepackage{soul}
\usepackage{tikz}
\usepackage{amsthm}
\usepackage{amsmath}
\usepackage{amssymb}
\usepackage{mathrsfs}
\usepackage{geometry}
\usepackage{graphicx}
\usepackage{amsfonts}
\usepackage{epstopdf}
\usepackage{placeins}
\usepackage{enumerate}
\usepackage{enumitem}
\usepackage{subcaption}
\usepackage{hyperref}
\usepackage[normalem]{ulem} % enables \sout, avoids changing \emph

\usetikzlibrary{calc}

\numberwithin{equation}{section}

\allowdisplaybreaks

\newcommand{\Rmnum}[1]{\uppercase\expandafter{\romannumeral#1}} % Uppercase roman number

\def\Xint#1{\mathchoice
{\XXint\displaystyle\textstyle{#1}}%
{\XXint\textstyle\scriptstyle{#1}}%
{\XXint\scriptstyle\scriptscriptstyle{#1}}%
{\XXint\scriptscriptstyle\scriptscriptstyle{#1}}%
\!\int}
\def\XXint#1#2#3{{\setbox0=\hbox{$#1{#2#3}{\int}$ }
\vcenter{\hbox{$#2#3$ }}\kern-.6\wd0}}

\def\dashint{\Xint-}

\newtheorem{theorem}{Theorem}[section]
\newtheorem{proposition}[theorem]{Proposition}
\newtheorem{lemma}[theorem]{Lemma}
\newtheorem{corollary}[theorem]{Corollary}
\newtheorem{remark}[theorem]{Remark}

\newtheorem{asmA}{Assumption}

\newtheorem{asmAp}{Assumption}

\renewcommand{\thefootnote}{}

\AtEndDocument{
\bigskip{\footnotesize%
   \textsc{Department of Mathematics, Aarhus University, 8000 Aarhus C, Denmark} \par
  \textit{E-mail address}: \texttt{yangmengqh@gmail.com}
}}

\begin{document}

%\title{\large{This is title}}
\title{On the dichotomy of $p$-walk dimensions on metric measure spaces}
\author{Meng Yang}
\date{}

\maketitle

\abstract{On a volume doubling metric measure space endowed with a family of $p$-energies such that the Poincar\'e inequality and the cutoff Sobolev inequality with $p$-walk dimension $\beta_p$ hold, for $p$ in an open interval $I\subseteq (1,+\infty)$, we prove the following dichotomy: either $\beta_p=p$ for \emph{all} $p\in I$, or $\beta_p>p$ for \emph{all} $p\in I$.}

\footnote{\textsl{Date}: \today}
\footnote{The author is grateful to Fabrice Baudoin for stimulating discussions. The author appreciates helpful anonymous comments on a previous version of this work, in particular on the relative compactness of metric balls.}

\renewcommand{\thefootnote}{\arabic{footnote}}
\setcounter{footnote}{0}

\section{Introduction}

%The paper \cite{GHL15}, \cite{Shi24a} will be used.

On many fractals, including the Sierpi\'nski gasket and the Sierpi\'nski carpet, there exists a diffusion with a heat kernel satisfying the following two-sided sub-Gaussian estimates:
\begin{align*}\label{eq_HKbeta}\tag*{HK($\beta$)}
\frac{C_1}{V(x,t^{1/\beta})}\exp\left(-C_2\left(\frac{d(x,y)}{t^{1/\beta}}\right)^{\frac{\beta}{\beta-1}}\right)\le p_t(x,y)\le\frac{C_3}{V(x,t^{1/\beta})}\exp\left(-C_4\left(\frac{d(x,y)}{t^{1/\beta}}\right)^{\frac{\beta}{\beta-1}}\right),\nonumber
\end{align*}
where $\beta$ is a new parameter called the walk dimension, which is always strictly greater than 2 on fractals. For example, $\beta=\frac{\log5}{\log2}$ on the Sierpi\'nski gasket (see \cite{BP88,Kig89}), $\beta\approx2.09697$ on the Sierpi\'nski carpet (see \cite{BB89,BB90,BBS90,BB92,KZ92,HKKZ00}). For $\beta=2$, \ref{eq_HKbeta} is indeed the classical Gaussian estimates.

By the standard Dirichlet form theory, a diffusion corresponds to a local regular Dirichlet form (see \cite{FOT11}). The Dirichlet form framework generalizes the classical Dirichlet integral $\int_{\mathbb{R}^d}|\nabla f(x)|^2 \mathrm{d} x$ in $\mathbb{R}^d$. For general $p>1$, extending the classical $p$-energy $\int_{\mathbb{R}^d}|\nabla f(x)|^p \mathrm{d} x$ in $\mathbb{R}^d$, as initiated by \cite{HPS04}, the study of $p$-energy on fractals and general metric measure spaces has been recently advanced considerably, see \cite{CGQ22,Shi24,BC23,MS23,Kig23,CGYZ24,AB25,AES25a}. In this setting, a new parameter $\beta_p$, called the $p$-walk dimension, naturally arises in connection with a $p$-energy. Notably, $\beta_2$ coincides with $\beta$ in \ref{eq_HKbeta}.

Since $\beta_2$ is typically strictly greater than $2$ on many classical fractals, it is natural to expect that $\beta_p$ would be strictly greater than $p$ on these fractals as well. On the Vicsek set, $\beta_p=p+d_h-1>p$, where $d_h=\frac{\log5}{\log3}$ is the Hausdorff dimension; see \cite{BC23}. On the Sierpi\'nski gasket and the Sierpi\'nski carpet, the inequality $\beta_p>p$ was established in \cite{KS24a}, whereas the exact value of $\beta_p$ remains unknown, except for $\beta_2=\frac{\log5}{\log2}$ on the Sierpi\'nski gasket. The main motivation of this paper is to study the behavior of the inequality $\beta_p>p$ in a more systematic way. More precisely, under the volume doubling condition, assume that the Poincar\'e inequality and the cutoff Sobolev inequality with $p$-walk dimension $\beta_p$ hold for all $p$ in an open interval $I\subseteq (1,+\infty)$. We prove that either $\beta_p=p$ for all $p\in I$, or $\beta_p>p$ for all $p\in I$; see Theorem \ref{thm_dichotomy}. Consequently, if $2\in I$ or $I=(1,+\infty)$---which is usually the case---the inequality $\beta_2>2$ suffices to obtain the corresponding strict inequality for all $p\in I$.

We provide a brief outline of the proof as follows. Firstly, under the volume doubling condition, the Poincar\'e inequality and the capacity upper bound with $p$-walk dimension $\beta_p$, the quotient $\alpha_p=\frac{\beta_p}{p}$ can be characterized in terms of the critical exponent of certain Besov spaces, see \cite{Shi24a}. Utilizing this characterization, we obtain regularity properties of the functions $p\mapsto \alpha_p$ and $p\mapsto\beta_p$. In particular, $\alpha_p\ge1$ is monotone decreasing and continuous in $p$, while $\beta_p$ is monotone increasing and continuous in $p$, see \cite{Bau24}. This implies that $\beta_p\ge p$ for all $p$, and that the set $\{p:\beta_p=p\}=\{p:\alpha_p=1\}$ is a relatively closed subinterval of $I$ of the form $[p,+\infty)\cap I$. Secondly, assume that $\{p:\alpha_p=1\}$ is non-empty. Take any $p$ in this set, then $\beta_p=p$. By adapting the techniques in \cite{KM20} to the $p$-energy setting, we prove that the conjunction of the Poincar\'e inequality and the cutoff Sobolev inequality with $p$-walk dimension $\beta_p=p$ implies that the associated $p$-energy measure is absolutely continuous with respect to the underlying measure, and that the associated intrinsic metric is bi-Lipschtiz equivalent to the underlying metric, see Theorem \ref{thm_AC}. In this case, by adapting the techniques in \cite{Stu94,Stu95,KSZ12,KZ12} to the $p$-energy setting, we obtain that Lipschitz functions are ``locally" contained in the domain of the $p$-energy, see Theorem \ref{thm_Lip}, and that a certain $(1,p)$-Poincar\'e inequality \hyperlink{eq_PI_1p}{$\mathrm{PI}_{\mathrm{Lip}}(1,p)$} holds. A very deep result from \cite{KZ08} further provides that such $(1,p)$-Poincar\'e inequality is an open ended condition, hence there exists $\varepsilon>0$ such that \hyperlink{eq_PI_1p}{$\mathrm{PI}_{\mathrm{Lip}}(1,q)$} holds for any $q>p-\varepsilon$, which in turn implies that the critical exponent $\alpha_q=1$ for any $q>p-\varepsilon$. Therefore, $\{p:\alpha_p=1\}$ is open in $I$. In summary, $\{p:\alpha_p=1\}$ is both relatively open and relatively closed in $I$; hence the dichotomy follows directly.

Throughout this paper, the letters $C$, $C_1$, $C_2$, $C_A$, $C_B$ will always refer to some positive constants and may change at each occurrence. The sign $\asymp$ means that the ratio of the two sides is bounded from above and below by positive constants. The sign $\lesssim$ ($\gtrsim$) means that the LHS is bounded by positive constant times the RHS from above (below). We use $x_+$ to denote the positive part of $x\in \mathbb{R}$, that is, $x_+=\max\{x,0\}$. For two $\sigma$-finite Borel measures $\mu$, $\nu$, the notion $\mu\le\nu$ means that $\mu\ll \nu$ and $\frac{\mathrm{d}\mu}{\mathrm{d}\nu}\le1$, that is $\mu$ is absolutely continuous with respect to $\nu$ with Radon-Nikodym derivative bounded by 1. We use $\#A$ to denote the cardinality of a set $A$.

%We will use the notation $\lfloor x\rfloor$ ($\lceil x\rceil$) to denote the largest integer less than or equal to (the smallest integer greater than or equal to) $x\in \mathbb{R}$.

\section{Statement of main results}

Let $(X,d,m)$ be a \emph{complete} metric measure space, that is, $(X,d)$ is a complete locally compact separable metric space and $m$ is a positive Radon measure on $X$ with full support. Throughout this paper, we always assume that all metric balls are relatively compact. For any $x\in X$, for any $r\in(0,+\infty)$, denote $B(x,r)=\{y\in X:d(x,y)<r\}$ and $V(x,r)=m(B(x,r))$. If $B=B(x,r)$, then denote $\delta B=B(x,\delta r)$ for any $\delta\in(0,+\infty)$. Let $\mathcal{B}(X)$ be the family of all Borel measurable subsets of $X$. Let $C(X)$ be the family of all continuous functions on $X$. Let $C_c(X)$ be the family of all continuous functions on $X$ with compact support. Denote $\dashint_A=\frac{1}{m(A)}\int_A$ and $u_A=\dashint_Au \mathrm{d} m$ for any measurable set $A$ with $m(A)\in(0,+\infty)$ and any function $u$ such that the integral $\int_Au \mathrm{d}m$ is well-defined.

Let $\varepsilon\in(0,+\infty)$. We say that $V$ is an $\varepsilon$-net (of $(X,d)$) if $V\subseteq X$ satisfies that for any distinct $x,y\in V$, we have $d(x,y)\ge \varepsilon$, and for any $z\in X$, there exists $x\in V$ such that $d(x,z)< \varepsilon$. Since $(X,d)$ is separable, all $\varepsilon$-nets are countable.

We say that the chain condition \ref{eq_CC} holds if there exists $C_{cc}\in(0,+\infty)$ such that for any $x,y\in X$, for any positive integer $n$, there exists a sequence $\{x_k:0\le k\le n\}$ of points in $X$ with $x_0=x$ and $x_n=y$ such that
\begin{equation*}\label{eq_CC}\tag*{CC}
d(x_k,x_{k-1})\le C_{cc} \frac{d(x,y)}{n}\text{ for any }k=1,\ldots,n.
\end{equation*}
Throughout this paper, we always assume \ref{eq_CC}.

We say that the volume doubling condition \ref{eq_VD} holds if there exists $C_{VD}\in(0,+\infty)$ such that
\begin{equation*}\label{eq_VD}\tag*{VD}
V(x,2r)\le C_{VD}V(x,r)\text{ for any }x\in X,r\in(0,+\infty).
\end{equation*}

%We say that the volume regular condition \ref{eq_VPhi} holds if there exists $C_{VR}\in(0,+\infty)$ such that
%\begin{equation*}\label{eq_VPhi}\tag*{V($\Phi$)}
%\frac{1}{C_{VR}}\Phi(r)\le V(x,r)\le C_{VR}\Phi(r)\text{ for any }x\in X,r\in(0,+\infty).
%\end{equation*}
%For $d_h\in(0,+\infty)$, we say that the Ahlfors regular condition \hypertarget{eq_Vdh}{V($d_h$)} holds if \ref{eq_VPhi} holds with $\Phi:r\mapsto r^{d_h}$.

We say that $(\mathcal{E},\mathcal{F})$ is a $p$-energy on $(X,d,m)$ if $\mathcal{F}$ is a dense subspace of $L^p(X;m)$ and $\mathcal{E}:\mathcal{F}\to[0,+\infty) $ satisfies the following conditions.

\begin{enumerate}[label=(\arabic*)]
\item $\mathcal{E}^{1/p}$ is a semi-norm on $\mathcal{F}$, that is, for any $f,g\in\mathcal{F}$, $c\in\mathbb{R}$, we have $\mathcal{E}(f)\ge0$, $\mathcal{E}(cf)^{1/p}=|c|\mathcal{E}(f)^{1/p}$ and $\mathcal{E}(f+g)^{1/p}\le\mathcal{E}(f)^{1/p}+\mathcal{E}(g)^{1/p}$.
\item (Closed property) $(\mathcal{F},\mathcal{E}(\cdot)^{1/p}+\lVert {\cdot}\rVert_{L^p(X;m)})$ is a Banach space.
\item (Markovian property) For any $\varphi\in C(\mathbb{R})$ with $\varphi(0)=0$ and $|\varphi(t)-\varphi(s)|\le|t-s|$ for any $t$, $s\in\mathbb{R}$, for any $f\in\mathcal{F}$, we have $\varphi(f)\in\mathcal{F}$ and $\mathcal{E}(\varphi(f))\le\mathcal{E}(f)$.
\item (Regular property) $\mathcal{F}\cap C_c(X)$ is uniformly dense in $C_c(X)$ and $(\mathcal{E}(\cdot)^{1/p}+\lVert {\cdot}\rVert_{L^p(X;m)})$-dense in $\mathcal{F}$.
\item (Strongly local property) For any $f,g\in\mathcal{F}$ with compact support and $g$ constant in an open neighborhood of $\mathrm{supp}(f)$, we have $\mathcal{E}(f+g)=\mathcal{E}(f)+\mathcal{E}(g)$.
\item ($p$-Clarkson's inequality) For any $f,g\in\mathcal{F}$, we have
\begin{equation*}\label{eq_Cla}\tag*{Cla}
\begin{cases}
\mathcal{E}(f+g)+\mathcal{E}(f-g)\ge2 \left(\mathcal{E}(f)^{\frac{1}{p-1}}+\mathcal{E}(g)^{\frac{1}{p-1}}\right)^{p-1}&\text{if }p\in(1,2],\\
\mathcal{E}(f+g)+\mathcal{E}(f-g)\le2 \left(\mathcal{E}(f)^{\frac{1}{p-1}}+\mathcal{E}(g)^{\frac{1}{p-1}}\right)^{p-1}&\text{if }p\in[2,+\infty).\\
\end{cases}
\end{equation*}
\end{enumerate}
Moreover, we also always assume the following condition.
\begin{itemize}
\item ($\mathcal{F}\cap L^\infty(X;m)$ is an algebra) For any $f,g\in\mathcal{F}\cap L^\infty(X;m)$, we have $fg\in\mathcal{F}$ and
\begin{equation*}\label{eq_Alg}\tag*{Alg}
\mathcal{E}(fg)^{{1}/{p}}\le \lVert {f}\rVert_{L^\infty(X;m)}\mathcal{E}(g)^{1/p}+\lVert {g}\rVert_{L^\infty(X;m)}\mathcal{E}(f)^{1/p}.
\end{equation*}
\end{itemize}
Denote $\mathcal{E}_\lambda(\cdot)=\mathcal{E}(\cdot)+\lambda\lVert\cdot \rVert^p_{L^p(X;m)}$ for any $\lambda\in(0,+\infty)$. Indeed, a general condition called the generalized $p$-contraction property was introduced in \cite{KS24a}, which implies \ref{eq_Cla}, \ref{eq_Alg}, and holds on a large family of metric measure spaces.

By \cite[Theorem 1.4]{Sas25}, a $p$-energy $(\mathcal{E},\mathcal{F})$ corresponds to a (canonical) $p$-energy measure $\Gamma:\mathcal{F}\times\mathcal{B}(X)\to[0,+\infty)$, $(f,A)\mapsto\Gamma(f)(A)$ satisfying the following conditions.
\begin{enumerate}[label=(\arabic*)]
\item For any $f\in\mathcal{F}$, $\Gamma(f)(\cdot)$ is a positive Radon measure on $X$ with $\Gamma(f)(X)=\mathcal{E}(f)$.
\item For any $A\in\mathcal{B}(X)$, $\Gamma(\cdot)(A)^{1/p}$ is a semi-norm on $\mathcal{F}$.
\item For any $f,g\in\mathcal{F}\cap C_c(X)$, $A\in\mathcal{B}(X)$, if $f-g$ is constant on $A$, then $\Gamma(f)(A)=\Gamma(g)(A)$.
\item ($p$-Clarkson's inequality) For any $f,g\in\mathcal{F}$, for any $A\in\mathcal{B}(X)$, we have
\begin{equation*}
\begin{cases}
\Gamma(f+g)(A)+\Gamma(f-g)(A)\ge2 \left(\Gamma(f)(A)^{\frac{1}{p-1}}+\Gamma(g)(A)^{\frac{1}{p-1}}\right)^{p-1}&\text{if }p\in(1,2],\\
\Gamma(f+g)(A)+\Gamma(f-g)(A)\le2 \left(\Gamma(f)(A)^{\frac{1}{p-1}}+\Gamma(g)(A)^{\frac{1}{p-1}}\right)^{p-1}&\text{if }p\in[2,+\infty).\\
\end{cases}
\end{equation*}
\item (Chain rule) For any $f\in\mathcal{F}\cap C_c(X)$, for any piecewise $C^1$ function $\varphi:\mathbb{R}\to\mathbb{R}$, we have $\mathrm{d}\Gamma(\varphi(f))=|\varphi'(f)|^p\mathrm{d}\Gamma(f)$.
\end{enumerate}
Using the chain rule, we have the following condition.
\begin{itemize}
\item (Strong sub-additivity) For any $f,g\in\mathcal{F}$, we have $f\vee g$, $f\wedge g\in\mathcal{F}$ and
\begin{equation*}\label{eq_SubAdd}\tag*{SubAdd}
\mathcal{E}(f\vee g)+\mathcal{E}(f\wedge g)\le\mathcal{E}(f)+\mathcal{E}(g).
\end{equation*}
\end{itemize}

Let
\begin{align*}
&\mathcal{F}_{\mathrm{loc}}=\left\{u:
\begin{array}{l}
\text{for any relatively compact open set }U,\\
\text{there exists }u^\#\in \mathcal{F}\text{ such that }u=u^\#\text{ }m\text{-a.e. in }U
\end{array}
\right\}.
\end{align*}
For any $u\in \mathcal{F}_{\mathrm{loc}}$, let $\Gamma(u)|_U=\Gamma(u^\#)|_U$, where $u^\#$, $U$ are given as above, then $\Gamma(u)$ is a well-defined positive Radon measure on $X$. By the strongly local property of $(\mathcal{E},\mathcal{F})$, we have the following result:
\begin{equation}\label{eq_local}
\text{If }u,v\in \mathcal{F}_{\mathrm{loc}}\text{ satisfy that }\Gamma(u)\le m, \Gamma(v)\le m,\text{ then }\Gamma(u\vee v)\le m.
\end{equation}

Let $\Psi:[0,+\infty)\to[0,+\infty)$ be a doubling function, that is, $\Psi$ is a homeomorphism, which implies that $\Psi$ is strictly increasing continuous and $\Psi(0)=0$, and there exists $C_\Psi\in(1,+\infty)$, called a doubling constant of $\Psi$, such that $\Psi(2r)\le C_\Psi\Psi(r)$ for any $r\in(0,+\infty)$.

We say that the Poincar\'e inequality \ref{eq_PI} holds if there exist $C_{PI}\in(0,+\infty)$, $A_{PI}\in[1,+\infty)$ such that for any ball $B$ with radius $r\in(0,+\infty)$, for any $f\in \mathcal{F}$, we have
\begin{equation*}\label{eq_PI}\tag*{PI$_p$($\Psi$)}
\int_B\lvert f-f_B\rvert^p \mathrm{d} m\le C_{PI}\Psi(r)\int_{A_{PI}B} \mathrm{d}\Gamma(f).
\end{equation*}
For $\beta_p\in(0,+\infty)$, we say that the Poincar\'e inequality \hypertarget{eq_PIbeta}{PI$_p$($\beta_p$)} holds if \ref{eq_PI} holds with $\Psi:r\mapsto r^{\beta_p}$.

Let $U$, $V$ be two open subsets of $X$ satisfying $U\subseteq \overline{U}\subseteq V$. We say that $\phi\in \mathcal{F}$ is a cutoff function for $U\subseteq V$ if $0\le\phi\le1$ in $X$, $\phi=1$ in an open neighborhood of $\overline{U}$ and $\mathrm{supp}(\phi)\subseteq V$, where $\mathrm{supp}(f)$ refers to the support of the measure of $|f|\mathrm{d} m$ for any given function $f$.

We say that the cutoff Sobolev inequality \ref{eq_CS} holds if there exist $C_{1}$, $C_{2}\in(0,+\infty)$, $A_{S}\in(1,+\infty)$ such that for any ball $B(x,r)$, there exists a cutoff function $\phi\in \mathcal{F}$ for $B(x,r)\subseteq B(x,A_Sr)$ such that for any $f\in \mathcal{F}$, we have
\begin{equation*}\label{eq_CS}\tag*{CS$_p$($\Psi$)}
\int_{B(x,A_{S}r)}|\widetilde{f}|^p \mathrm{d}\Gamma(\phi)\le C_{1}\int_{B(x,A_{S}r)}\mathrm{d}\Gamma(f)+\frac{C_{2}}{\Psi(r)}\int_{B(x,A_{S}r)}|f|^p \mathrm{d}m,
\end{equation*}
where $\widetilde{f}$ is a quasi-continuous modification of $f$, such that $\widetilde{f}$ is uniquely determined $\Gamma(\phi)$-a.e. in $X$, see \cite[Section 8]{Yan25a} for more details. For $\beta_p\in(0,+\infty)$, we say that the cutoff Sobolev inequality \hypertarget{eq_CSbeta}{CS$_p$($\beta_p$)} holds if \ref{eq_CS} holds with $\Psi:r\mapsto r^{\beta_p}$.

Let $A_1,A_2\in \mathcal{B}(X)$. We define the capacity between $A_1$, $A_2$ as
\begin{align*}
&\mathrm{cap}(A_1,A_2)=\inf\left\{\mathcal{E}(\varphi):\varphi\in \mathcal{F},
\begin{array}{l}
\varphi=1\text{ in an open neighborhood of }A_1,\\
\varphi=0\text{ in an open neighborhood of }A_2
\end{array}
\right\},
\end{align*}
here we use the convention that $\inf\emptyset=+\infty$.

We say that the capacity upper bound \ref{eq_ucap} holds if there exist $C_{cap}\in(0,+\infty)$, $A_{cap}\in(1,+\infty)$ such that for any ball $B(x,r)$, we have
\begin{align*}
\mathrm{cap}\left(B(x,r),X\backslash B(x,A_{cap}r)\right)&\le C_{cap} \frac{V(x,r)}{\Psi(r)}.\label{eq_ucap}\tag*{$\text{cap}_p(\Psi)_{\le}$}
%\mathrm{cap}\left(B(x,r),X\backslash B(x,A_{cap}r)\right)&\ge \frac{1}{C_{cap}} \frac{V(x,r)}{\Psi(r)}.\label{eq_lcap}\tag*{$\text{cap}(\Psi)_{\ge}$}
\end{align*}
For $\beta_p\in(0,+\infty)$, we say that \hypertarget{eq_ucapbeta}{$\text{cap}_p(\beta_p)_{\le}$}  holds if \ref{eq_ucap} holds with $\Psi:r\mapsto r^{\beta_p}$. Under \ref{eq_VD}, by taking $f\equiv1$ in $B(x,A_Sr)$, it is easy to see that \ref{eq_CS} (resp. \hyperlink{eq_CSbeta}{$\text{CS}_p(\beta_p)$}) implies \ref{eq_ucap} (resp. \hyperlink{eq_ucapbeta}{$\text{cap}_p(\beta_p)_{\le}$}).

The main result of this paper is the following dichotomy.

\begin{theorem}\label{thm_dichotomy}
Assume \ref{eq_VD}. Let $I\subseteq(1,+\infty)$ be an open interval. Assume for any $p\in I$, there exists a $p$-energy $(\mathcal{E},\mathcal{F})$ such that \hyperlink{eq_PIbeta}{PI$_p$($\beta_p$)}, \hyperlink{eq_CSbeta}{CS$_p$($\beta_p$)} hold. Then
\begin{enumerate}[label=(\roman*)]
\item either $\beta_p=p$ for all $p\in I$,
\item or $\beta_p>p$ for all $p\in I$.
\end{enumerate}
\end{theorem}

As a direct corollary, we obtain the strict inequality $\beta_p>p$ for all $p\in(1,+\infty)$ on the Sierpi\'nski gasket and the Sierpi\'nski carpet as follows.

\begin{corollary}
On the Sierpi\'nski gasket and the Sierpi\'nski carpet, for any $p\in(1,+\infty)$, let $(\mathcal{E},\mathcal{F})$ be the $p$-energy with $p$-walk dimension $\beta_p$, as constructed in \cite{HPS04,CGQ22} for the Sierpi\'nski gasket, and in \cite{Shi24,MS23} for the Sierpi\'nski carpet. Then $\beta_p>p$ for any $p\in(1,+\infty)$.
\end{corollary}

\begin{proof}
For any $p\in(1,+\infty)$, by \cite[Corollary 2.5]{Yan25c}, \hyperlink{eq_PIbeta}{PI$_p$($\beta_p$)}, \hyperlink{eq_CSbeta}{CS$_p$($\beta_p$)} hold on the Sierpi\'nski gasket; by \cite[Corollary 2.10]{Yan25d}, \hyperlink{eq_PIbeta}{PI$_p$($\beta_p$)}, \hyperlink{eq_CSbeta}{CS$_p$($\beta_p$)} hold on the Sierpi\'nski carpet. By the standard and widely known result that $\beta_2>2$ on these fractals, see for instance \cite{BP88,BB89,Kig89}, the result follows.
\end{proof}

\begin{remark}
This result was also obtained in \cite[Theorem 9.8 and Theorem 9.13]{KS24a}, where the proof relies on the self-similar property. The contribution of our work is that once \hyperlink{eq_PIbeta}{PI$_p$($\beta_p$)}, \hyperlink{eq_CSbeta}{CS$_p$($\beta_p$)} are established---which is the case on many fractals and metric measure spaces, see \cite[Theorem 2.3]{Yan25c} and \cite[Theorem 2.9]{Yan25d} for several equivalent characteriza-tions---the proof of $\beta_p>p$ for all $p$ could be reduced to proving $\beta_2>2$, which would be much easier to handle. Indeed, such an argument can be applied to a family of strongly symmetric p.c.f. self-similar sets, and to a family of $p$-conductively homogeneous compact metric spaces, see \cite[Remark 2.6]{Yan25c} and the references therein.
\end{remark}

Let us introduce the key ingredients for the proof. The intrinsic metric $\rho:X\times X\to[0,+\infty]$ of $(\mathcal{E},\mathcal{F})$ is given by
\begin{equation}\label{eq_intrinsic}
\rho(x,y)=\sup \left\{f(x)-f(y):f\in \mathcal{F}_{\mathrm{loc}}\cap C(X),\Gamma(f)\le m\right\}.
\end{equation}
By definition, $\rho$ is only a pseudo metric and not necessarily a metric. However, under the following assumption:

\begin{asmAp}\label{asm_Ap}
The topology induced by $\rho$ is equivalent to the original topology on $(X,d)$.
\end{asmAp}

\noindent
we have $\rho$ is indeed a metric, as a consequence of the remark after \cite[Assumption (A')]{Stu95} and the fact that $X$ is connected, which in turn follows from \ref{eq_CC} and \cite[PROPOSITION A.1]{KM20}. We will also need another stronger assumption as follows:
\begin{asmA}\label{asm_A}
$\rho$ is a complete metric on $X$ which is compatible with the original topology on $(X,d)$.
\end{asmA}

\noindent
Assuming \ref{asm_A}, the metric balls with respect to $\rho$ are relatively compact; this property will be crucial in the proof of Proposition \ref{prop_geodesic} and the subsequent results. For a comparison between \ref{asm_A} and \ref{asm_Ap}, see \cite[Theorem 2]{Stu95}.

The first ingredient for the proof is that under \ref{asm_A}, Lipschitz functions with respect to $\rho$ are contained in $\mathcal{F}_{\mathrm{loc}}$. This result parallels \cite[Theorem 2.1]{KZ12} in the Dirichlet form setting.

We now introduce the related notions with respect to $\rho$. Let $B_\rho(x,r)=\{y\in X:\rho(x,y)<r\}$ be the open ball centered at $x$ of radius $r$ with respect to $\rho$. For $x\in X$, for a function $u$ defined in an open neighborhood of $x$, its pointwise Lipschitz constant at $x$ with respect to $\rho$ is defined as
$$\mathrm{Lip}_\rho u(x)=\lim_{r\downarrow0}\sup_{y:\rho(x,y)\in(0,r)}\frac{\lvert u(x)-u(y)\rvert}{\rho(x,y)}.$$
Let $V$ be an open subset of $X$. We say a function $u$ defined in $V$ is Lipschitz in $V$ with respect to $\rho$ if there exists $K\in(0,+\infty)$ such that $\lvert u(x)-u(y)\rvert\le K\rho(x,y)$ for any $x$, $y\in V$. Let $\mathrm{Lip}_\rho(V)$ be the family of all Lipschitz functions in $V$ with respect to $\rho$ and
$$\lVert u\rVert_{\mathrm{Lip}_\rho(V)}=\sup_{x,y\in V,x\ne y}\frac{\lvert u(x)-u(y)\rvert}{\rho(x,y)}\text{ for any }u\in \mathrm{Lip}_\rho(V).$$

\begin{theorem}\label{thm_Lip}
Assume \ref{asm_A} and that $(X,\rho,m)$ satisfies the volume doubling condition, that is, there exists $C\in(0,+\infty)$ such that
\begin{equation}\label{eq_rho_VD}
m(B_\rho(x,2r))\le Cm(B_\rho(x,r))\text{ for any }x\in X, r\in(0,+\infty).
\end{equation}
Then $\mathrm{Lip}_\rho(X)\subseteq \mathcal{F}_{\mathrm{loc}}$ and $\Gamma(u)\le (\mathrm{Lip}_\rho u)^p m$ for any $u\in \mathrm{Lip}_\rho (X)$.
\end{theorem}

The second ingredient for the proof is the absolute continuity of the $p$-energy measure with respect to the underlying measure, and the bi-Lipschitz equivalence between the intrinsic metric and the underlying metric, as stated below.

A $\sigma$-finite Borel measure $\mu$ on $X$ is called a minimal energy-dominant measure of $(\mathcal{E},\mathcal{F})$ if the following two conditions are satisfied.
\begin{enumerate}[label=(\roman*)]
\item (Domination) For any $f\in \mathcal{F}$, we have $\Gamma(f)\ll \mu$.
\item (Minimality) If another $\sigma$-finite Borel measure $\nu$ on $X$ also satisfies the above domination condition, then $\mu\ll \nu$.
\end{enumerate}
See \cite[Lemma 9.20]{MS23} for the existence of such a measure, and also \cite[Lemma 2.2]{Nak85}, \cite[LEMMAS 2.2, 2.3 and 2.4]{Hin10} for the existence in the Dirichlet form setting.

\begin{theorem}\label{thm_AC}
Assume \ref{eq_VD}, \ref{eq_PI}, \ref{eq_CS} and
\begin{equation}\label{eq_Psi_AC}
\varlimsup_{r\downarrow0}\frac{\Psi(r)}{r^p}>0.
\end{equation}
Then $m$ is a minimal energy-dominant measure of $(\mathcal{E},\mathcal{F})$, hence $\Gamma(f)\ll m$ for any $f\in \mathcal{F}$. Moreover, $\rho$ is a geodesic metric on $X$, and $\rho$ is bi-Lipschitz equivalent to $d$, that is, there exists $C\in(0,+\infty)$ such that
$$\frac{1}{C}d(x,y)\le \rho(x,y)\le Cd(x,y)\text{ for any }x,y\in X.$$
In particular, assume \ref{eq_VD}, \hyperlink{eq_PIbeta}{PI$_p$($p$)}, \hyperlink{eq_CSbeta}{CS$_p$($p$)}, then all the above results hold.
\end{theorem}

\begin{remark}
We will follow an argument from \cite{KM20}, where the case $p=2$ was considered.
\end{remark}

This paper is organized as follows. In Section \ref{sec_intrinsic}, we prove Theorem \ref{thm_Lip}. In Section \ref{sec_AC}, we prove Theorem \ref{thm_AC}. In Section \ref{sec_dichotomy}, we prove Theorem \ref{thm_dichotomy}.

\section{Proof of Theorem \ref{thm_Lip}}\label{sec_intrinsic}

Let $(\mathcal{E},\mathcal{F})$ be a $p$-energy with intrinsic metric $\rho$ given as in Equation (\ref{eq_intrinsic}). Let $\rho(x,\cdot):y\mapsto\rho(x,y)$ be the distance function to $x$ with respect to $\rho$.

Firstly, we present the following two results in the $p$-energy setting, which are parallel to \cite[Lemma 1']{Stu94} and \cite[Lemma 3, Theorem 1]{Stu95} in the Dirichlet form setting, respectively. These results show that, under \ref{asm_A}, the distance functions $\rho(x,\cdot)$ belong to $\mathcal{F}_{\mathrm{loc}}$, and that $\rho$ is a geodesic metric.

\begin{proposition}\label{prop_rhox}
Assume \ref{asm_Ap}. For any $x\in X$, the distance function $\rho(x,\cdot):y\mapsto \rho(x,y)$ satisfies that $\rho(x,\cdot)\in \mathcal{F}_{\mathrm{loc}}\cap C(X)$ and $\Gamma(\rho(x,\cdot))\le m$.
\end{proposition}

\begin{proof}
By assumption, we have $(X,\rho)$ is separable, for any $n\ge1$, let $\{z^{(n)}_i\}_{i\ge1}$ be a $\frac{1}{n}$-net of $(X,\rho)$. For any $i\ge1$, by definition, there exists $\psi^{(n)}_i\in \mathcal{F}_{\mathrm{loc}}\cap C(X)$ with $\Gamma(\psi^{(n)}_i)\le m$ such that
\begin{equation}\label{eq_rhox1}
\rho(x,z^{(n)}_i)-\frac{1}{n}<\psi^{(n)}_i(x)-\psi^{(n)}_i(z^{(n)}_i)\le\rho(x,z^{(n)}_i).
\end{equation}
Moreover, for any $y\in B_\rho(z^{(n)}_i,\frac{1}{n})$, we have
$$\frac{1}{n}>\rho(y,z^{(n)}_i)\ge \psi^{(n)}_i(y)-\psi^{(n)}_i(z^{(n)}_i),$$
which gives
$$\psi^{(n)}_i(y)\le\psi^{(n)}_i(z^{(n)}_i)+\frac{1}{n}\overset{\text{Eq. (\ref{eq_rhox1})}}{\scalebox{6}[1]{$<$}}\psi^{(n)}_i(x)-\rho(x,z^{(n)}_i)+\frac{2}{n}<\psi^{(n)}_i(x)-\rho(x,y)+\frac{3}{n},$$
hence $\psi^{(n)}_i(x)-\psi^{(n)}_i\ge\rho(x,\cdot)-\frac{3}{n}$ in $B_\rho(z^{(n)}_i,\frac{1}{n})$. Since $\psi^{(n)}_i(x)-\psi^{(n)}_i\le \rho(x,\cdot)$ in $X$, let $\phi^{(n)}_i=(\psi^{(n)}_i(x)-\psi^{(n)}_i)_+$, then
\begin{equation}\label{eq_rhox2}
\phi^{(n)}_i\in \mathcal{F}_{\mathrm{loc}}\cap C(X)\text{ and }\Gamma(\phi^{(n)}_i)\le m,
\end{equation}
\begin{equation}\label{eq_rhox3}
0\le \phi^{(n)}_i\le \rho(x,\cdot)\text{ in }X,
\end{equation}
\begin{equation}\label{eq_rhox4}
\phi^{(n)}_i\ge \rho(x,\cdot)-\frac{3}{n}\text{ in }B_\rho(z^{(n)}_i,\frac{1}{n}).
\end{equation}
By replacing $\phi^{(n)}_i$ with $\max_{1\le j\le i}\phi^{(n)}_j$, we may assume that $\phi^{(n)}_i$ is increasing in $i$. By Equation (\ref{eq_local}), $\phi^{(n)}_i$ satisfies Equation (\ref{eq_rhox2}), moreover, $\phi^{(n)}_i$ satisfies Equation (\ref{eq_rhox3}) and
\begin{equation}\label{eq_rhox5}
\phi^{(n)}_j\ge \rho(x,\cdot)-\frac{3}{n}\text{ in }B_\rho(z^{(n)}_i,\frac{1}{n})\text{ for any }j\ge i\ge 1.
\end{equation}

For any relatively compact open subset $X_0\subseteq X$, by \ref{asm_Ap}, there exists $M>0$ such that $X_0\subseteq \overline{X_0}\subseteq B_\rho(x,M)$. By the regular property of $(\mathcal{E},\mathcal{F})$, there exists $\psi\in \mathcal{F}\cap C_c(X)$ with $0\le \psi\le 1$ in $X$, $\psi=1$ in $B_\rho(x,M)$, and $\mathrm{supp}(\psi)\subseteq B_\rho(x,2M)$. Let $\varphi^{(n)}_i=\phi^{(n)}_i\wedge (M\psi)$, then $\varphi^{(n)}_i\in \mathcal{F}\cap C_c(X)$, $\varphi^{(n)}_i=\phi^{(n)}_i$ in $B_\rho(x,M)$, and $\mathrm{supp}(\varphi^{(n)}_i)\subseteq B_\rho(x,2M)$. It is obvious that $\{\varphi_i^{(n)}\}_{i\ge1}$ is $L^p(X;m)$-bounded. Since
\begin{align*}
&\mathcal{E}(\varphi^{(n)}_i)=\Gamma(\varphi^{(n)}_i)(B_\rho(x,2M))\\
&\le \Gamma(\phi^{(n)}_i)(B_\rho(x,2M))+\Gamma(M\psi)(B_\rho(x,2M))\\
&\overset{\text{Eq. (\ref{eq_rhox2})}}{\scalebox{6}[1]{$\le$}}m(B_\rho(x,2M))+M^p \mathcal{E}(\psi),
\end{align*}
we have $\{\varphi_i^{(n)}\}_{i\ge1}$ is $\mathcal{E}$-bounded, which gives $\{\varphi_i^{(n)}\}_{i\ge1}$ is $\mathcal{E}_1$-bounded. By the Banach–Alao-glu theorem (see \cite[Theorem 3 in Chapter 12]{Lax02}), there exists a subsequence, still denoted by $\{\varphi^{(n)}_i\}_{i\ge1}$, which is $\mathcal{E}_1$-weakly-convergent to some element $\phi^{(n)}\in \mathcal{F}$. By Mazur's lemma, here we refer to the version in \cite[Theorem 2 in Section V.1]{Yos95}, for any $i\ge1$, there exist $I_i\ge i$, $\lambda_k^{(i)}\ge0$ for $k=i,\ldots,I_i$ with $\sum_{k=i}^{I_i}\lambda_k^{(i)}=1$ such that $\{\sum_{k=i}^{I_i}\lambda_k^{(i)}\varphi^{(n)}_k\}_{i\ge1}$ is $\mathcal{E}_1$-convergent to $\phi^{(n)}$. For any $i\ge1$, by Equation (\ref{eq_rhox3}), we have
$$0\le\sum_{k=i}^{I_i}\lambda_k^{(i)}\varphi^{(n)}_k=\sum_{k=i}^{I_i}\lambda_k^{(i)}\phi^{(n)}_k\le \rho(x,\cdot)\text{ in }B_\rho(x,M),$$
hence $0\le\phi^{(n)}\le\rho(x,\cdot)$ in $B_\rho(x,M)$; moreover, for any $j\ge i\ge1$, by Equation (\ref{eq_rhox5}), we have
$$\sum_{k=j}^{I_j}\lambda_k^{(j)}\varphi^{(n)}_k=\sum_{k=j}^{I_j}\lambda_k^{(j)}\phi^{(n)}_k\ge \rho(x,\cdot)-\frac{3}{n}\text{ in }B_\rho(z^{(n)}_i,\frac{1}{n})\cap B_\rho(x,M),$$
hence $\phi^{(n)}\ge\rho(x,\cdot)-\frac{3}{n}$ in $B_{\rho}(z^{(n)}_i,\frac{1}{n})\cap B_\rho(x,M)$ for any $i\ge1$, which gives $\phi^{(n)}\ge\rho(x,\cdot)-\frac{3}{n}$ in $B_\rho(x,M)$. Since $\varphi^{(n)}_i=\phi^{(n)}_i$ in $B_\rho(x,M)$, by Equation (\ref{eq_rhox2}), we have $\Gamma(\varphi^{(n)}_i)\le m$ in $B_\rho(x,M)$, by the triangle inequality for $\Gamma(\cdot)(A)^{1/p}$ for any $A\in \mathcal{B}(X)$, we have $\Gamma(\sum_{k=i}^{I_i}\lambda_k^{(i)}\varphi^{(n)}_k)\le m$ in $B_\rho(x,M)$, which gives $\Gamma(\phi^{(n)})\le m$ in $B_\rho(x,M)$.

Hence, for any $n\ge1$, there exists $\phi^{(n)}\in \mathcal{F}$ satisfying that $\rho(x,\cdot)-\frac{3}{n}\le\phi^{(n)}\le\rho(x,\cdot)$ in $B_\rho(x,M)$, and $\Gamma(\phi^{(n)})\le m$ in $B_\rho(x,M)$. Similar to the above argument, let $\eta\in \mathcal{F}\cap C_c(X)$ satisfy $0\le\eta\le 1$ in $X$, $\eta=1$ on $\overline{X_0}$, and $\mathrm{supp}(\eta)\subseteq B_\rho(x,M)$, then certain convex combinations of $\{\phi^{(n)}\wedge (M\eta)\}_{n\ge1}$ is $\mathcal{E}_1$-convergent to some $\phi\in \mathcal{F}$, where $\Gamma(\phi)\le m$ in $X_0$ and $\phi=\rho(x,\cdot)$ in $X_0$. Therefore, $\rho(x,\cdot)\in \mathcal{F}_{\mathrm{loc}}\cap C(X)$ satisfies $\Gamma(\rho(x,\cdot))\le m$.
\end{proof}

\begin{proposition}\label{prop_geodesic}
Assume \ref{asm_A}. For any $x$, $y\in X$, let $R=\rho(x,y)<+\infty$, for any $r\in [0,R]$, there exists $z\in X$ such that $\rho(x,z)=r$, $\rho(z,y)=R-r$. Hence $(X,\rho)$ is a geodesic space.
\end{proposition}

\begin{proof}
Without loss of generality, we may assume that $R=\rho(x,y)\in(0,+\infty)$, $r\in(0,R)$. Suppose there exist $x$, $y$, $r$ such that no such $z$ exists, then the closed balls $\overline{B_\rho(x,r)}$, $\overline{B_\rho(y,R-r)}$ are disjoint. By \cite[Theorem 2]{Stu95}, assuming \ref{asm_A}, $\overline{B_\rho(x,r)}$, $\overline{B_\rho(y,R-r)}$ are compact, hence with respect to $\rho$, their distance $D=\mathrm{dist}_\rho(\overline{B_\rho(x,r)},\overline{B_\rho(y,R-r)})\in(0,+\infty)$. Let $\delta\in(0,\frac{1}{3}D)$, then $B_\rho(x,r+\delta)\cap B_\rho(y,R-r+\delta)=\emptyset$, let
$$f=
\begin{cases}
(r+\delta)-\rho(x,\cdot)&\text{in }B_\rho(x,r+\delta),\\
\rho(y,\cdot)-(R-r+\delta)&\text{in }B_\rho(y,R-r+\delta),\\
0&\text{otherwise}.
\end{cases}
$$
Then by Proposition \ref{prop_rhox}, we have $f\in \mathcal{F}_{\mathrm{loc}}\cap C(X)$, and by the strongly local property of $(\mathcal{E},\mathcal{F})$, we have $\Gamma(f)=1_{B_\rho(x,r+\delta)}\Gamma(\rho(x,\cdot))+1_{B_\rho(y,R-r+\delta)}\Gamma(\rho(y,\cdot))\le m$, hence
$$\rho(x,y)\ge f(x)-f(y)=(r+\delta)+(R-r+\delta)=R+2\delta>R=\rho(x,y),$$
contradiction. In particular, for any $x$, $y\in X$, there exists $z\in X$ such that $\rho(x,z)=\rho(z,y)=\frac{1}{2}\rho(x,y)$. By \ref{asm_A}, $(X,\rho)$ is complete, hence $(X,\rho)$ is a geodesic space, see for instance \cite[Remarks 1.4 (1)]{BH99}.
\end{proof}

Secondly, we present the following two preparatory results for the proof of Theorem \ref{thm_Lip}.

\begin{lemma}[{\cite[LEMMA 6.30]{Che99}, \cite[Lemma 2.3]{KZ12}}]\label{lem_Che}
Assume \ref{asm_A} and that $(X,\rho,m)$ satisfies the volume doubling condition Equation (\ref{eq_rho_VD}). Then for any ball $B_\rho(x_0,r_0)$, there exists $C\in[1,+\infty)$ such that for any $n\ge1$, for any $u\in \mathrm{Lip}_\rho(B_\rho(x_0,r_0))$, there exists a finite family of mutually disjoint balls $\{B_\rho(x_{n,i},r_{n,i})\}_i$ with $x_{n,i}\in B_\rho(x_0,r_0)$ and $r_{n,i}\le r_0$ for any $i$, such that
\begin{equation}\label{eq_Che1}
\mathrm{dist}_\rho \left(B_\rho(x_{n,i},r_{n,i}),B_\rho(x_{n,j},r_{n,j})\right)\ge \frac{1}{2}\left(r_{n,i}+r_{n,j}\right)\text{ for any }i\ne j,
\end{equation}
\begin{equation}\label{eq_Che2}
m \left(B_\rho(x_{0},r_0)\backslash \cup_i B_\rho(x_{n,i},r_{n,i})\right)\le \frac{C}{n^p} m \left(B_\rho(x_0,r_0)\right),
\end{equation}
\begin{equation}\label{eq_Che3}
\left(\frac{1}{m(B_\rho(x_{n,i},3r_{n,i}))}\int_{B_\rho(x_{n,i},3r_{n,i})}\lvert \mathrm{Lip}_\rho u-\mathrm{Lip}_\rho u(x_{n,i})\rvert^p \mathrm{d}m\right)^{1/p}\le \frac{1}{n},
\end{equation}
\begin{equation}\label{eq_Che4}
\frac{\lvert u(x)-u(y)\rvert}{\rho(x,y)}\le \mathrm{Lip}_\rho u(x_{n,i})+\frac{1}{n}\text{ for any }x,y\in B_\rho(x_{n,i},r_{n,i})\text{ with }\rho(x,y)\ge\frac{1}{n}r_{n,i}.
\end{equation}
\end{lemma}

\begin{remark}
Assuming \ref{asm_A}, Proposition \ref{prop_geodesic} gives that $(X,\rho)$ is a geodesic space; hence \cite[LEMMA 6.30]{Che99} applies.
\end{remark}

We need the following result to extend a Lipschitz function from a subset to the whole space. This result parallels \cite[Lemma 2.2]{KZ12} in the Dirichlet form setting.

\begin{lemma}\label{lem_VtoX}
Assume \ref{asm_A}. Let $V$ be a bounded open subset of $(X,\rho)$. For any $v\in \mathrm{Lip}_\rho(V)$ with $\lVert v\rVert_{\mathrm{Lip}_\rho(V)}\le 1$, let
$$u=\sup_{z\in V}\left\{v(z)-\rho(z,\cdot)\right\}.$$
Then $u=v$ in $V$, $u\in \mathcal{F}_{\mathrm{loc}}\cap \mathrm{Lip}_\rho(X)$, $\lVert u\rVert_{\mathrm{Lip}_\rho(X)}\le 1$ and $\Gamma(u)\le m$.
\end{lemma}

\begin{proof}
It is obvious that $u=v$ in $V$, $u\in \mathrm{Lip}_\rho(X)$, and $\lVert u\rVert_{\mathrm{Lip}_\rho(X)}\le1$. Let $D=\mathrm{diam}_\rho(V)<+\infty$. For any $n\ge1$, let $\{z_{n,i}\}_i$ be a $(\frac{1}{n}D)$-net of $(V,\rho)$, which is a finite set. Let
$$u_n=\max_i \left\{v(z_{n,i})-\rho(z_{n,i},\cdot)\right\}.$$
Then for any $x\in X$, by definition, we have $u_n(x)\le u(x)$. For any $z\in V$, there exists $z_{n,i}$ such that $\rho(z,z_{n,i})<\frac{1}{n}D$, hence
\begin{align*}
&u_n(x)\ge v(z_{n,i})-\rho(z_{n,i},x)\\
&\ge \left(v(z)-\rho(z,z_{n,i})\right)-\left(\rho(z,x)+\rho(z,z_{n,i})\right)\\
&\ge \left(v(z)-\rho(z,x)\right)-\frac{2}{n}D,
\end{align*}
taking the supremum with respect to $z$, we have $u_n(x)\ge u(x)-\frac{2}{n}D$ for any $x\in X$. Therefore, $\{u_n\}_n$ converges uniformly to $u$.

By Proposition \ref{prop_rhox} and Equation (\ref{eq_local}), we have $u_n\in \mathcal{F}_{\mathrm{loc}}$ and $\Gamma(u_n)\le m$. For any bounded open subset $U$ of $(X,\rho)$, we have $\{\int_{U}\mathrm{d}\Gamma(u_n)\}_n$ is bounded, and $\{u_n\}_n$ is $L^p(U;m)$-convergent to $u$. Let
$$\mathcal{F}^{\mathrm{ref}}(U)=\left\{u\in \mathcal{F}_{\mathrm{loc}}\cap L^p(U;m):\int_U \mathrm{d}\Gamma(u)<+\infty\right\},$$
then $(\mathcal{F}^{\mathrm{ref}}(U),\mathcal{E}_1)$ is a reflexive Banach space. Since $\{u_n\}_n$ is a bounded sequence in $(\mathcal{F}^{\mathrm{ref}}(U),\mathcal{E}_1)$, by the Banach–Alaoglu theorem (see \cite[Theorem 3 in Chapter 12]{Lax02}), there exists a subsequence, still denoted by $\{u_n\}_{n}$, which is $\mathcal{E}_1$-weakly-convergent to some element $w\in\mathcal{F}^{\mathrm{ref}}(U)$. By Mazur's lemma, here we refer to the version in \cite[Theorem 2 in Section V.1]{Yos95}, for any $n\ge1$, there exist $I_n\ge n$, $\lambda_k^{(n)}\ge0$ for $k=n,\ldots,I_n$ with $\sum_{k=n}^{I_n}\lambda_k^{(n)}=1$ such that $\{\sum_{k=n}^{I_n}\lambda_k^{(n)}u_k\}_{n}$ is $\mathcal{E}_1$-convergent to $w$. Since
$$\lVert\sum_{k=n}^{I_n}\lambda_k^{(n)}u_k-u\rVert_{L^p(U;m)}\le\sum_{k=n}^{I_n}\lambda_k^{(n)}\lVert u_k-u\rVert_{L^p(U;m)}\le\sup_{k\ge n}\lVert u_k-u\rVert_{L^p(U;m)}\to0$$
as $n\to+\infty$, we have $u=w$ in $U$, $u\in \mathcal{F}_{\mathrm{loc}}$. By the triangle inequality, we have
\begin{align*}
&\left(\int_U \mathrm{d}\Gamma(u)\right)^{1/p}=\left(\int_U \mathrm{d}\Gamma(w)\right)^{1/p}=\lim_{n\to+\infty}\left(\int_U \mathrm{d}\Gamma \left(\sum_{k=n}^{I_n}\lambda_k^{(n)}u_k\right)\right)^{1/p}\\
&\le\varliminf_{n\to+\infty}\sum_{k=n}^{I_n}\lambda_k^{(n)}\left(\int_U \mathrm{d}\Gamma(u_k)\right)^{1/p}\le \varliminf_{n\to+\infty}\sum_{k=n}^{I_n}\lambda_k^{(n)}m(U)^{1/p}=m(U)^{1/p},
\end{align*}
hence $\Gamma(u)(U)\le m(U)$ for any bounded open subset $U$, which gives $\Gamma(u)\le m$.
\end{proof}

We give the proof of Theorem \ref{thm_Lip} as follows.

\begin{proof}[Proof of Theorem \ref{thm_Lip}]
Our argument follows the MacShane extension technique, as in the proof of \cite[Theorem 2.1]{KZ12}. Let $u\in \mathrm{Lip}_\rho(X)$ and $L=\lVert u\rVert_{\mathrm{Lip}_\rho(X)}$. We only need to show that for any ball $B_\rho(x_0,r_0)$, there exists $v\in \mathcal{F}_{\mathrm{loc}}$ such that $v=u$ in $B_\rho(x_0,r_0)$, and
$$\int_{B_\rho(x_0,r_0)}\mathrm{d}\Gamma(v)\le\int_{B_\rho(x_0,r_0)}(\mathrm{Lip}_\rho u)^p \mathrm{d}m.$$

For any $n\ge1$, let $u_n\in \mathrm{Lip}_\rho(X)$ be given as follows. Let $\{B_\rho(x_{n,i},r_{n,i})\}_i$ be given as in Lemma \ref{lem_Che}, and $L_{n,i}=\mathrm{Lip}_\rho u(x_{n,i})+\frac{1}{n}\le L+\frac{1}{n}$. For any $i$, let $\{z_{n,i,k}\}_{k}$ be a $\left(\frac{1}{n}r_{n,i}\right)$-net of $(B_\rho(x_{n,i},r_{n,i}),\rho)$, which is a finite set. Let
$$u_n=\max_k \left\{u(z_{n,i,k})-L_{n,i}\rho(z_{n,i,k},\cdot)\right\}\text{ in }B_\rho(x_{n,i},r_{n,i}),$$
then it is obvious that $\lVert u_n\rVert_{\mathrm{Lip}_\rho(B_\rho(x_{n,i},r_{n,i}))}\le L_{n,i}$. For any $x\in B_\rho(x_{n,i},r_{n,i})$, by definition, there exists $k$ such that $u_n(x)=u(z_{n,i,k})-L_{n,i}\rho(z_{n,i,k},x)$. Since $(X,\rho)$ is a geodesic space, there exists a geodesic $\gamma$ connecting $x$ and $z_{n,i,k}$, then
\begin{align*}
&\mathrm{Lip}_\rho u_n(x)\ge\varlimsup_{\gamma\ni y\to x}\frac{u_n(y)-u_n(x)}{\rho(x,y)}\\
&\ge \varlimsup_{\gamma\ni y\to x}\frac{\left(u(z_{n,i,k})-L_{n,i}\rho(z_{n,i,k},y)\right)-\left(u(z_{n,i,k})-L_{n,i}\rho(z_{n,i,k},x)\right)}{\rho(x,y)}\\
&=\varlimsup_{\gamma\ni y\to x}L_{n,i}\frac{\rho(z_{n,i,k},x)-\rho(z_{n,i,k},y)}{\rho(x,y)}\overset{\gamma:\text{ geodesic}}{\scalebox{6}[1]{$=$}}L_{n,i}.
\end{align*}
Hence $\mathrm{Lip}_\rho u_n\equiv\lVert u_n\rVert_{\mathrm{Lip}_\rho(B_\rho(x_{n,i},r_{n,i}))}=L_{n,i}$ in $B_\rho(x_{n,i},r_{n,i})$. By Proposition \ref{prop_rhox} and Equation (\ref{eq_local}), we have $\Gamma(u_n)\le L_{n,i}^pm=(\mathrm{Lip}_\rho u_n)^pm$ in $B_\rho(x_{n,i},r_{n,i})$, hence
$$\Gamma(u_n)\le (\mathrm{Lip}_\rho u_n)^pm\text{ in }\bigcup_iB_\rho(x_{n,i},r_{n,i}).$$
By Equation (\ref{eq_Che4}), we have $\frac{\lvert u(z_{n,i,k})-u(z_{n,i,l})\rvert}{\rho(z_{n,i,k},z_{n,i,l})}\le L_{n,i}$ for any $k\ne l$, which gives that $u_n(z_{n,i,k})=u(z_{n,i,k})$.

We claim that $\lVert u_n\rVert_{\mathrm{Lip}_\rho(\cup_iB_\rho(x_{n,i},r_{n,i}))}\le7L+2$. Indeed, for any $x\in B_\rho(x_{n,i},r_{n,i})$, $y\in B_\rho(x_{n,j},r_{n,j})$ with $x\ne y$, if $i=j$, then $\frac{\lvert u_n(x)-u_n(y)\rvert}{\rho(x,y)}\le \lVert u_n\rVert_{\mathrm{Lip}_\rho(B_\rho(x_{n,i},r_{n,i}))}=L_{n,i}\le L+1$; if $i\ne j$, then by Equation (\ref{eq_Che1}), we have
$$\rho(x,y)\ge \mathrm{dist}_\rho \left(B_\rho(x_{n,i},r_{n,i}),B_\rho(x_{n,j},r_{n,j})\right)\ge \frac{1}{2}\left(r_{n,i}+r_{n,j}\right).$$
There exist $k$, $l$ such that $\rho(x,z_{n,i,k})<\frac{1}{n}r_{n,i}\le r_{n,i}$, $\rho(y,z_{n,j,l})< \frac{1}{n}r_{n,j}\le r_{n,j}$. Since
$$\rho(z_{n,i,k},z_{n,j,l})\le \mathrm{dist}_\rho \left(B_\rho(x_{n,i},r_{n,i}),B_\rho(x_{n,j},r_{n,j})\right)+2(r_{n,i}+r_{n,j})\le 5\rho(x,y),$$
and $u_n(z_{n,i,k})=u(z_{n,i,k})$, $u_n(z_{n,j,l})=u(z_{n,j,l})$, we have
\begin{align*}
&\lvert u_n(x)-u_n(y)\rvert\\
&\le \lvert u_n(x)-u_n(z_{n,i,k})\rvert+\lvert u(z_{n,i,k})-u(z_{n,j,l})\rvert+\lvert u_n(y)-u_n(z_{n,j,l})\rvert\\
&\le L_{n,i}r_{n,i}+L\rho(z_{n,i,k},z_{n,j,l})+L_{n,j}r_{n,j}\\
&\le (L+1)(r_{n,i}+r_{n,j})+5L\rho(x,y)\\
&\le (7L+2)\rho(x,y).
\end{align*}
Hence $\lVert u_n\rVert_{\mathrm{Lip}_\rho(\cup_iB_\rho(x_{n,i},r_{n,i}))}\le7L+2$.

Let
$$u_n=\sup_{z\in \cup_{i}B_\rho(x_{n,i},r_{n,i})}\left\{u_n(z)-\lVert u_n\rVert_{\mathrm{Lip}_\rho(\cup_iB_\rho(x_{n,i},r_{n,i}))}\rho(z,\cdot)\right\},$$
then $u_n\in \mathrm{Lip}_\rho(X)$ is well-defined and
$$\lVert u_n\rVert_{\mathrm{Lip}_\rho(X)}=\lVert u_n\rVert_{\mathrm{Lip}_\rho(\cup_iB_\rho(x_{n,i},r_{n,i}))}\le 7L+2.$$
By Lemma \ref{lem_VtoX}, we have $u_n\in \mathcal{F}_{\mathrm{loc}}$ and $\Gamma(u_n)\le (7L+2)^pm$, which gives $\{\int_{B_\rho(x_0,r_0)}\mathrm{d}\Gamma(u_n)\}_n$  is bounded.

We claim that $\{u_n\}_{n}$ is $L^p(B_\rho(x_0,r_0);m)$-convergent to $u$. Indeed, for arbitrary $x\in B_\rho(x_{n,i},r_{n,i})$, there exists $z_{n,i,k}$ such that $\rho(x,z_{n,i,k})<\frac{1}{n}r_{n,i}$, recall that $u_n(z_{n,i,k})=u(z_{n,i,k})$, hence
\begin{align*}
&\lvert u_n(x)-u(x)\rvert\\
&\le \lvert u_n(x)-u_n(z_{n,i,k})\rvert+\lvert u(x)-u(z_{n,i,k})\rvert\\
&\le \left(L_{n,i}+L\right)\frac{1}{n}r_{n,i}\le \frac{(2L+1)r_0}{n},
\end{align*}
which gives
$$|u_n(x)-u(x)|\le \frac{(2L+1)r_0}{n}\text{ for any }x\in \cup_{i}B_\rho(x_{n,i},r_{n,i}).$$
For any $x\in B_\rho(x_0,r_0)\backslash\cup_iB_\rho(x_{n,i},r_{n,i})$, take arbitrary $B_\rho(x_{n,i},r_{n,i})$ and arbitrary $z_{n,i,k}\in B_\rho(x_{n,i},r_{n,i})$, then
\begin{align*}
&\lvert u_n(x)-u(x)\rvert\\
&\le \lvert u_n(x)-u_n(z_{n,i,k})\rvert+\lvert u(x)-u(z_{n,i,k})\rvert\\
&\le\left(7L+2+L\right)\rho(x,z_{n,i,k})\overset{(*)}{\scalebox{3}[1]{$\le$}}3(8L+2)r_0,
\end{align*}
where in ($*$), we use the fact that $\rho(x,z_{n,i,k})\le\rho(x,x_{n,i})+\rho(x_{n,i},z_{n,i,k})\le 3r_0$. Hence
\begin{align*}
&\int_{B_\rho(x_0,r_0)}|u_n-u|^p \mathrm{d}m\\
&= \int_{B_\rho(x_0,r_0)\cap\cup_i B_\rho(x_{n,i},r_{n,i})}|u_n-u|^p \mathrm{d}m+\int_{B_\rho(x_0,r_0)\backslash\cup_i B_\rho(x_{n,i},r_{n,i})}|u_n-u|^p \mathrm{d}m\\
&\le \left(\frac{(2L+1)r_0}{n}\right)^p m(B_\rho(x_0,r_0))+(3(8L+2)r_0)^pm(B_\rho(x_0,r_0)\backslash\cup_i B_\rho(x_{n,i},r_{n,i}))\\
&\overset{\text{Eq. (\ref{eq_Che2})}}{\scalebox{6}[1]{$\le$}}\frac{1}{n^p}\left(3^p(8L+2)^pC_1+(2L+1)^p\right){r_0^pm(B_\rho(x_0,r_0))}\to0
\end{align*}
as $n\to+\infty$, where $C_1$ is the constant appearing in Equation (\ref{eq_Che2}), which gives $\{u_n\}_{n}$ is $L^p(B_\rho(x_0,r_0);m)$-convergent to $u$.

Let
$$\mathcal{F}^{\mathrm{ref}}(B_\rho(x_0,r_0))=\left\{u\in \mathcal{F}_{\mathrm{loc}}\cap L^p(B_\rho(x_0,r_0);m):\int_{B_\rho(x_0,r_0)}\mathrm{d}\Gamma(u)<+\infty\right\},$$
then $(\mathcal{F}^{\mathrm{ref}}(B_\rho(x_0,r_0)),\mathcal{E}_1)$ is a reflexive Banach space. Since $\{u_n\}_n$ is a bounded sequence in $(\mathcal{F}^{\mathrm{ref}}(B_\rho(x_0,r_0)),\mathcal{E}_1)$, by the Banach–Alaoglu theorem (see \cite[Theorem 3 in Chapter 12]{Lax02}), there exists a subsequence, still denoted by $\{u_n\}_{n}$, which is $\mathcal{E}_1$-weakly-convergent to some element $v\in\mathcal{F}^{\mathrm{ref}}(B_\rho(x_0,r_0))$. By Mazur's lemma, here we refer to the version in \cite[Theorem 2 in Section V.1]{Yos95}, for any $n\ge1$, there exist $I_n\ge n$, $\lambda_k^{(n)}\ge0$ for $k=n,\ldots,I_n$ with $\sum_{k=n}^{I_n}\lambda_k^{(n)}=1$ such that $\{\sum_{k=n}^{I_n}\lambda_k^{(n)}u_k\}_{n}$ is $\mathcal{E}_1$-convergent to $v$. Since
\begin{align*}
&\lVert\sum_{k=n}^{I_n}\lambda_k^{(n)}u_k-u\rVert_{L^p(B_\rho(x_0,r_0);m)}\\
&\le\sum_{k=n}^{I_n}\lambda_k^{(n)}\lVert u_k-u\rVert_{L^p(B_\rho(x_0,r_0);m)}\\
&\le\sup_{k\ge n}\lVert u_k-u\rVert_{L^p(B_\rho(x_0,r_0);m)}\to0
\end{align*}
as $n\to+\infty$, we have $u=v$ in $B_\rho(x_0,r_0)$. By the triangle inequality, we have
\begin{align*}
&\left(\int_{B_\rho(x_0,r_0)}\mathrm{d}\Gamma(v)\right)^{1/p}=\lim_{n\to+\infty}\left(\int_{B_\rho(x_0,r_0)}\mathrm{d}\Gamma \left(\sum_{k=n}^{I_n}\lambda_k^{(n)}u_k\right)\right)^{1/p}\\
&\le\varliminf_{n\to+\infty}\sum_{k=n}^{I_n}\lambda_k^{(n)}\left(\int_{B_\rho(x_0,r_0)}\mathrm{d}\Gamma \left(u_k\right)\right)^{1/p}\le \varlimsup_{n\to+\infty}\left(\int_{B_\rho(x_0,r_0)}\mathrm{d}\Gamma \left(u_n\right)\right)^{1/p}.
\end{align*}
For any $n$, we have
\begin{align*}
&\lvert\left(\int_{B_\rho(x_0,r_0)\cap \cup_iB_\rho(x_{n,i},r_{n,i})}(\mathrm{Lip}_\rho u)^p \mathrm{d}m\right)^{1/p}-\left(\int_{B_\rho(x_0,r_0)\cap\cup_iB_\rho(x_{n,i},r_{n,i})}(\mathrm{Lip}_\rho u_n)^p \mathrm{d}m\right)^{1/p}\rvert\\
&\le \left(\int_{B_\rho(x_0,r_0)\cap\cup_iB_\rho(x_{n,i},r_{n,i})}\lvert\mathrm{Lip}_\rho u-\mathrm{Lip}_\rho u_n\rvert^p \mathrm{d}m\right)^{1/p}\\
&\overset{\{B_\rho(x_{n,i},r_{n,i})\}_i:\text{ disjoint}}{\underset{\mathrm{Lip}_\rho u_n\equiv L_{n,i}\text{ in }B_\rho(x_{n,i},r_{n,i})}{\scalebox{16}[1]{=}}}\left(\sum_i\int_{B_\rho(x_0,r_0)\cap B_\rho(x_{n,i},r_{n,i})}\lvert\mathrm{Lip}_\rho u-L_{n,i}\rvert^p \mathrm{d}m\right)^{1/p}\\
&\overset{L_{n,i}=\mathrm{Lip}_\rho u(x_{n,i})+\frac{1}{n}}{\scalebox{12}[1]{$\le$}}\left(\sum_i\int_{B_\rho(x_{n,i},r_{n,i})}\lvert\mathrm{Lip}_\rho u-\mathrm{Lip}_\rho u(x_{n,i})\rvert^p \mathrm{d}m\right)^{1/p}\\
&\hspace{90pt}+\left(\int_{B_\rho(x_0,r_0)}\left(\frac{1}{n}\right)^p \mathrm{d}m\right)^{1/p}\\
&\overset{\text{Eq. (\ref{eq_Che3})}}{\scalebox{6}[1]{$\le$}}\left(\frac{1}{n^p}\sum_im(B_\rho(x_{n,i},3r_{n,i}))\right)^{1/p}+\frac{1}{n}m \left(B_\rho(x_0,r_0\right))^{1/p}\\
&\overset{\text{Eq. (\ref{eq_rho_VD})}}{\scalebox{6}[1]{$\le$}}\left(\frac{C_2^2}{n^p}\sum_im(B_\rho(x_{n,i},r_{n,i}))\right)^{1/p}+\frac{1}{n}m \left(B_\rho(x_0,r_0\right))^{1/p}\\
&\overset{\{B_\rho(x_{n,i},r_{n,i})\}_i:\text{ disjoint}}{\scalebox{15}[1]{$=$}}\frac{C_2^{2/p}}{n}m\left(\cup_iB_\rho(x_{n,i},r_{n,i})\right)^{1/p}+\frac{1}{n}m \left(B_\rho(x_0,r_0\right))^{1/p}\\
&\overset{x_{n,i}\in B(x_0,r_0)}{\underset{r_{n,i}\le r_0}{\scalebox{10}[1]{$\le$}}}\frac{C_2^{2/p}+1}{n}m \left(B_\rho(x_0,2r_0\right))^{1/p},
\end{align*}
where $C_2$ is the constant appearing in Equation (\ref{eq_rho_VD}), hence
\begin{align*}
&\int_{B_\rho(x_0,r_0)}\mathrm{d}\Gamma(u_n)\\
&\le\sum_i\int_{B_\rho(x_0,r_0)\cap B_\rho(x_{n,i},r_{n,i})}L_{n,i}^p\mathrm{d}m+(7L+2)^pm\left({B_\rho(x_0,r_0)\backslash \cup_iB_\rho(x_{n,i},r_{n,i})}\right)\\
&=\int_{B_\rho(x_0,r_0)\cap\cup_i B_\rho(x_{n,i},r_{n,i})}(\mathrm{Lip}_\rho u_n)^p\mathrm{d}m+(7L+2)^pm\left({B_\rho(x_0,r_0)\backslash \cup_iB_\rho(x_{n,i},r_{n,i})}\right)\\
&\le \left(\left(\int_{B_\rho(x_0,r_0)\cap \cup_iB_\rho(x_{n,i},r_{n,i})}(\mathrm{Lip}_\rho u)^p \mathrm{d}m\right)^{1/p}+\frac{C_2^{2/p}+1}{n}m \left(B_\rho(x_0,2r_0\right))^{1/p}\right)^{p}\\
&\hspace{10pt}+(7L+2)^pm\left({B_\rho(x_0,r_0)\backslash \cup_iB_\rho(x_{n,i},r_{n,i})}\right)\\
&\le \left(\left(\int_{B_\rho(x_0,r_0)}(\mathrm{Lip}_\rho u)^p \mathrm{d}m\right)^{1/p}+\frac{C_2^{2/p}+1}{n}m \left(B_\rho(x_0,2r_0\right))^{1/p}\right)^{1/p}\\
&\hspace{10pt}+(7L+2)^p \frac{C_1}{n^p}m(B_\rho(x_0,r_0)),
\end{align*}
where we use Equation (\ref{eq_Che2}) in the last inequality. Therefore,
$$\int_{B_\rho(x_0,r_0)}\mathrm{d}\Gamma(v)\le \varlimsup_{n\to+\infty}\int_{B_\rho(x_0,r_0)}\mathrm{d}\Gamma \left(u_n\right)\le \int_{B_\rho(x_0,r_0)}(\mathrm{Lip}_\rho u)^p \mathrm{d}m.$$
\end{proof}

\section{Proof of Theorem \ref{thm_AC}}\label{sec_AC}

We follow the argument given in \cite[Section 4]{KM20} in the Dirichlet form setting.

\begin{lemma}\label{lem_Psi_p}
Assume \ref{eq_VD}, \ref{eq_PI}, \ref{eq_ucap} and Equation (\ref{eq_Psi_AC}). Then there exist $r_1\in(0,\mathrm{diam}(X))$, $C\in(0,+\infty)$ such that
\begin{equation}\label{eq_Psi_p}
\frac{1}{C}r^p\le\Psi(r)\le Cr^p\text{ for any }r\in(0,r_1).
\end{equation}

\end{lemma}

\begin{proof}
By the proof of the lower bound in \cite[Proposition 2.1]{Yan25a}, there exists $C_1\in(0,+\infty)$ such that
\begin{equation}\label{eq_Psi_p1}
\frac{1}{C_1}\left(\frac{R}{r}\right)^p\le\frac{\Psi(R)}{\Psi(r)}\text{ for any }R,r\in(0,\mathrm{diam}(X))\text{ with }r\le R.
\end{equation}
By Equation (\ref{eq_Psi_AC}), there exist $C_2\in(0,+\infty)$, $\{r_n\}_{n\ge1}\subseteq(0,\mathrm{diam}(X))$ such that $r_n\downarrow0$ as $n\to+\infty$ and $\frac{\Psi(r_n)}{r_n^p}\ge \frac{1}{C_2}>0$ for any $n\ge1$.

For any $r\in(0,r_1)$, by Equation (\ref{eq_Psi_p1}), we have $\frac{\Psi(r)}{r^p}\le C_1 \frac{\Psi(r_1)}{r_1^p}$, and for any $n\ge1$ with $r_n<r$, we have $\frac{\Psi(r)}{r^p}\ge \frac{1}{C_1}\frac{\Psi(r_n)}{r_n^p}\ge \frac{1}{C_1C_2}$. Hence Equation (\ref{eq_Psi_p}) holds with $C=\max\{C_1C_2,C_1 \frac{\Psi(r_1)}{r_1^p}\}$.
\end{proof}

\begin{lemma}\label{lem_fxr}
Assume \ref{eq_VD}, \ref{eq_CS} and Equation (\ref{eq_Psi_p}). Then there exists $C\in(0,+\infty)$ such that for any $x\in X$, $r\in(0,r_1)$, let $f_{x,r}=(1-\frac{d(x,\cdot)}{r})_+$, then $f_{x,r}\in \mathcal{F}$ and $\Gamma(f_{x,r})\le\frac{C^p}{r^p} m$.
\end{lemma}

\begin{proof}
Let $C_1$ be the constant appearing in Equation (\ref{eq_Psi_p}). By \cite[Proposition 3.1]{Yan25c}, there exists $C_2\in(0,+\infty)$ such that for any $x\in X$, for any $r\in(0,r_1)$, for any $n\ge2$, for any $k=1,\ldots,n-1$, there exists a cutoff function $\phi_{n,k}\in \mathcal{F}$ for $B(x,\frac{k}{n}r)\subseteq B(x,\frac{k+1}{n}r)$ such that for any $g\in \mathcal{F}$, we have
\begin{align*}
&\int_{B(x,\frac{k+1}{n}r)\backslash \overline{B(x,\frac{k}{n}r)}}|\widetilde{g}|^p \mathrm{d}\Gamma(\phi_{n,k})\\
&\le \frac{1}{8}\int_{B(x,\frac{k+1}{n}r)\backslash \overline{B(x,\frac{k}{n}r)}}\mathrm{d}\Gamma(g)+\frac{C_2}{\Psi(\frac{1}{n}r)}\int_{B(x,\frac{k+1}{n}r)\backslash \overline{B(x,\frac{k}{n}r)}}|g|^p \mathrm{d}m\\
&\overset{\text{Eq. (\ref{eq_Psi_p})}}{\scalebox{6}[1]{$\le$}}\frac{1}{8}\int_{B(x,\frac{k+1}{n}r)\backslash \overline{B(x,\frac{k}{n}r)}}\mathrm{d}\Gamma(g)+\frac{C_1C_2n^p}{r^p}\int_{B(x,\frac{k+1}{n}r)\backslash \overline{B(x,\frac{k}{n}r)}}|g|^p \mathrm{d}m.
\end{align*}
For any $n\ge2$, let $\phi_n=\frac{1}{n-1}\sum_{k=1}^{n-1}\phi_{n,k}$, then $\phi_n\in \mathcal{F}$, $0\le\phi_n\le1$ in $X$, $\mathrm{supp}(\phi_n)\subseteq B(x,r)$, and $|\phi_n-f_{x,r}|\le \frac{2}{n}1_{B(x,r)}$ in $X$. By the strongly local property of $(\mathcal{E},\mathcal{F})$, for any $g\in \mathcal{F}$, we have
%$$\mathcal{E}(\phi_n)=\frac{1}{(n-1)^p}\sum_{k=1}^{n-1}\mathcal{E}(\phi_{n,k})\le\frac{1}{(n-1)^p}\frac{C_1C_2n^p}{r^p}V(x,r)\le$$
\begin{align*}
&\int_{B(x,r)}|\widetilde{g}|^p \mathrm{d}\Gamma(\phi_n)=\frac{1}{(n-1)^p}\sum_{k=1}^{n-1}\int_{B(x,\frac{k+1}{n}r)\backslash \overline{B(x,\frac{k}{n}r)}}|\widetilde{g}|^p \mathrm{d}\Gamma(\phi_{n,k})\\
&\le \frac{1}{(n-1)^p}\sum_{k=1}^{n-1}\left(\frac{1}{8}\int_{B(x,\frac{k+1}{n}r)\backslash \overline{B(x,\frac{k}{n}r)}}\mathrm{d}\Gamma(g)+\frac{C_1C_2n^p}{r^p}\int_{B(x,\frac{k+1}{n}r)\backslash \overline{B(x,\frac{k}{n}r)}}|g|^p \mathrm{d}m\right)\\
&\le\frac{1}{8(n-1)^p}\int_{B(x,r)}\mathrm{d}\Gamma(g)+\frac{C_1C_2n^p}{(n-1)^pr^p}\int_{B(x,r)}|g|^p \mathrm{d}m\\
&\le\frac{1}{8(n-1)^p}\int_{B(x,r)}\mathrm{d}\Gamma(g)+\frac{2^pC_1C_2}{r^p}\int_{B(x,r)}|g|^p \mathrm{d}m.
\end{align*}
By taking $g\equiv1$ in $B(x,r)$, we have $\mathcal{E}(\phi_n)\le\frac{2^pC_1C_2}{r^p}V(x,r)$. Since $\int_X|\phi_n|^p \mathrm{d}m\le V(x,r)$, we have $\{\phi_n\}_{n\ge2}$ is $\mathcal{E}_1$-bounded. Since $(\mathcal{F},\mathcal{E}_1)$ is a reflexive Banach space, by the Banach–Alaoglu theorem (see \cite[Theorem 3 in Chapter 12]{Lax02}), there exists a subsequence, still denoted by $\{\phi_n\}_{n\ge2}$, which is $\mathcal{E}_1$-weakly-convergent to some element $\phi\in \mathcal{F}$. By Mazur's lemma, here we refer to the version in \cite[Theorem 2 in Section V.1]{Yos95}, for any $n\ge2$, there exist $I_n\ge n$, $\lambda_k^{(n)}\ge0$ for $k=n,\ldots,I_n$ with $\sum_{k=n}^{I_n}\lambda_k^{(n)}=1$ such that $\{\sum_{k=n}^{I_n}\lambda_k^{(n)}\phi_k\}_{n\ge2}$ is $\mathcal{E}_1$-convergent to $\phi$. Since
$$\lvert\sum_{k=n}^{I_n}\lambda_k^{(n)}\phi_k-f_{x,r}\rvert\le\sum_{k=n}^{I_n}\lambda_k^{(n)}\lvert\phi_k-f_{x,r}\rvert\le\sum_{k=n}^{I_n}\lambda_k^{(n)}\frac{2}{k}1_{B(x,r)}\le \frac{2}{n}1_{B(x,r)}\to0$$
as $n\to+\infty$, we have $\{\sum_{k=n}^{I_n}\lambda_k^{(n)}\phi_k\}_{n\ge2}$ is $L^p(X;m)$-convergent to $f_{x,r}$, which gives $f_{x,r}=\phi\in \mathcal{F}$. For any $g\in \mathcal{F}\cap C_c(X)$, we have
\begin{align*}
&\left(\int_{B(x,r)}|{g}|^p \mathrm{d}\Gamma(f_{x,r})\right)^{1/p}\\
&\overset{(*)}{\scalebox{3}[1]{$=$}}\lim_{n\to+\infty}\left(\int_{B(x,r)}|{g}|^p \mathrm{d}\Gamma(\sum_{k=n}^{I_n}\lambda_k^{(n)}\phi_k)\right)^{1/p}\\
&\overset{(**)}{\scalebox{3}[1]{$\le$}}\varliminf_{n\to+\infty}\sum_{k=n}^{I_n}\lambda_k^{(n)}\left(\int_{B(x,r)}|{g}|^p \mathrm{d}\Gamma(\phi_k)\right)^{1/p}\le \varlimsup_{n\to+\infty}\left(\int_{B(x,r)}|{g}|^p \mathrm{d}\Gamma(\phi_n)\right)^{1/p}\\
&\le\varlimsup_{n\to+\infty}\left(\frac{1}{8(n-1)^p}\int_{B(x,r)}\mathrm{d}\Gamma(g)+\frac{2^pC_1C_2}{r^p}\int_{B(x,r)}|g|^p \mathrm{d}m\right)^{1/p}\\
&=\left(\frac{2^pC_1C_2}{r^p}\int_{B(x,r)}|g|^p \mathrm{d}m\right)^{1/p},
\end{align*}
where in ($*$), we use the fact that $g\in \mathcal{F}\cap C_c(X)$ is bounded, and in ($**$), we use the triangle inequality for $\left(\int_{B(x,r)}|g|^p \mathrm{d}\Gamma(\cdot)\right)^{1/p}$, hence
$$\int_{B(x,r)}|{g}|^p \mathrm{d}\Gamma(f_{x,r})\le\frac{2^pC_1C_2}{r^p}\int_{B(x,r)}|g|^p \mathrm{d}m\text{ for any }g\in \mathcal{F}\cap C_c(X).$$
By the regular property of $(\mathcal{E},\mathcal{F})$, we have $\Gamma(f_{x,r})\le \frac{C^p}{r^p}m$, where $C=2(C_1C_2)^{1/p}$.
\end{proof}

\begin{lemma}[Lipschitz partition of unity]\label{lem_partition}
Assume \ref{eq_VD}, \ref{eq_CS} and Equation (\ref{eq_Psi_p}). Then there exists $C\in(0,+\infty)$ such that for any $\varepsilon\in(0,\frac{r_1}{2})$, for any $\varepsilon$-net $V$, there exists a family of functions $\{\psi_z:z\in V\}\subseteq \mathcal{F}\cap C_c(X)$ satisfying the following conditions.
\begin{enumerate}[label=(CO\arabic*)]
\item\label{item_COunit} $\sum_{z\in V}\psi_z=1$.
\item\label{item_COspt} For any $z\in V$, $0\le\psi_z\le 1$ in $X$, and $\psi_z=0$ on $X\backslash B(z,2\varepsilon)$.
\item\label{item_COLip} For any $z\in V$, $\psi_z$ is $\frac{C}{\varepsilon}$-Lipschitz, that is, $\lvert \psi_z(x)-\psi_z(y)\rvert\le \frac{C}{\varepsilon}d(x,y)$ for any $x$, $y\in X$.
\item\label{item_COAC} For any $z\in V$, $\Gamma(\psi_z)\le \frac{C^p}{\varepsilon^p}m$.
\item\label{item_COenergy} For any $z\in V$, $\mathcal{E}(\psi_z)\le C\frac{V(z,\varepsilon)}{\varepsilon^p}$.
\end{enumerate}
\end{lemma}

\begin{proof}
Let $C_1$ be the constant appearing in Lemma \ref{lem_fxr}. By \ref{eq_VD}, there exists some positive integer $N$ depending only on $C_{VD}$ such that
$$\#\{z\in V:d(x,z)<4 \varepsilon\}\le N\text{ for any }x\in X.$$

For any $\varepsilon\in(0,\frac{r_1}{2})$, for any $\varepsilon$-net $V$, for any $z\in V$, let $f_{z,2 \varepsilon}\in \mathcal{F}$ be the function given by Lemma \ref{lem_fxr}. Then for any $x\in X$, there exists $z\in V$ such that $d(x,z)<\varepsilon$, hence $\sum_{z\in V}f_{z,2 \varepsilon}(x)\ge f_{z,2 \varepsilon}(x)\ge \frac{1}{2}$, and for any $z\in V$, if $f_{z,2 \varepsilon}(x)>0$, then $d(x,z)<2 \varepsilon$, hence $\sum_{z\in V}f_{z,2 \varepsilon}(x)=\sum_{z\in V:d(x,z)<2 \varepsilon}f_{z,2 \varepsilon}(x)\le \#\{z\in V:d(x,z)<2 \varepsilon\}\le N$. Therefore,
\begin{equation}\label{eq_sumf}
\frac{1}{2}\le\sum_{z\in V}f_{z,2 \varepsilon}\le N\text{ in }X.
\end{equation}

For any $z\in V$, let $\psi_z=\frac{f_{z,2 \varepsilon}}{\sum_{z\in V}f_{z,2 \varepsilon}}$, then $\psi_z\in C_c(X)$ is well-defined. It is obvious that \ref{item_COunit}, \ref{item_COspt} hold. By \cite[Proposition 2.3 (c)]{Shi24a}, we have $\psi_z\in \mathcal{F}$ and there exists some positive constant $C_2$ depending only on $p$, $N$ such that
\begin{align*}
&\mathcal{E}(\psi_z)=\Gamma(\psi_z)(B(z,2 \varepsilon))=\Gamma \left(\frac{f_{z,2 \varepsilon}}{\sum_{w\in V:d(z,w)<4 \varepsilon}f_{w,2 \varepsilon}}\right)(B(z,2 \varepsilon))\\
&\le C_2\sum_{w\in V:d(z,w)<4 \varepsilon}\Gamma(f_{w,2 \varepsilon})(B(z,2 \varepsilon))\overset{\text{Lem. \ref{lem_fxr}}}{\scalebox{6}[1]{$\le$}}C_2\sum_{w\in V:d(z,w)<4 \varepsilon}\frac{C_1^p}{(2 \varepsilon)^p}V(w,2 \varepsilon)\\
&\overset{\text{\ref{eq_VD}}}{\scalebox{3}[1]{$\le$}}\frac{C_1^pC_2C_{VD}^3N}{2^p}\frac{V(z,\varepsilon)}{\varepsilon^p},
\end{align*}
that is, \ref{item_COenergy} holds. Similarly, for any $z\in V$, for any $x\in X$, for any $r\in(0,2 \varepsilon)$, if $d(x,z)\ge 4 \varepsilon$, then $\Gamma(\psi_z)(B(x,r))=0$; if $d(x,z)<4 \varepsilon$, then
\begin{align*}
&\Gamma(\psi_z)(B(x,r))=\Gamma \left(\frac{f_{z,2 \varepsilon}}{\sum_{w\in V:d(x,w)<4 \varepsilon}f_{w,2 \varepsilon}}\right)(B(x,r))\le C_2\sum_{w\in V:d(x,w)<4 \varepsilon}\Gamma(f_{w,2 \varepsilon})(B(x,r))\\
&\overset{\text{Lem. \ref{lem_fxr}}}{\scalebox{6}[1]{$\le$}}C_2\sum_{w\in V:d(x,w)<4 \varepsilon}\frac{C_1^p}{(2 \varepsilon)^p}V(x,r)\le \frac{C_1^pC_2N}{2^p\varepsilon^p}V(x,r).
\end{align*}
Hence $\Gamma(\psi_z)\le \frac{C_1^pC_2N}{2^p \varepsilon^p}m$, that is, \ref{item_COAC} holds.

For any $z\in V$, for any $x,y\in X$, if $d(x,y)\ge 2\varepsilon$, then
$$|\psi_z(x)-\psi_z(y)|\le1\le \frac{1}{2\varepsilon}d(x,y).$$
If $d(x,y)<2\varepsilon$, recall that $|f_{w,2 \varepsilon}(x)-f_{w,2 \varepsilon}(y)|\le \frac{1}{2 \varepsilon}d(x,y)$ for any $w\in V$, then
\begin{align*}
&|\psi_z(x)-\psi_z(y)|=\lvert\frac{f_{z,2 \varepsilon}(x)}{\sum_{w\in V}f_{w,2 \varepsilon}(x)}-\frac{f_{z,2 \varepsilon}(y)}{\sum_{w\in V}f_{w,2 \varepsilon}(y)}\rvert\\
&\le \lvert\frac{f_{z,2 \varepsilon}(x)}{\sum_{w\in V}f_{w,2 \varepsilon}(x)}-\frac{f_{z,2 \varepsilon}(y)}{\sum_{w\in V}f_{w,2 \varepsilon}(x)}\rvert+\lvert\frac{f_{z,2 \varepsilon}(y)}{\sum_{w\in V}f_{w,2 \varepsilon}(x)}-\frac{f_{z,2 \varepsilon}(y)}{\sum_{w\in V}f_{w,2 \varepsilon}(y)}\rvert\\
&\le \frac{1}{\sum_{w\in V}f_{w,2 \varepsilon}(x)}\lvert f_{z,2 \varepsilon}(x)-f_{z,2 \varepsilon}(y)\rvert\\
&\hspace{10pt}+f_{z,2 \varepsilon}(y)\frac{1}{\sum_{w\in V}f_{w,2 \varepsilon}(x)}\frac{1}{\sum_{w\in V}f_{w,2 \varepsilon}(y)}\lvert\sum_{w\in V}f_{w,2 \varepsilon}(x)-\sum_{w\in V}f_{w,2 \varepsilon}(y)\rvert\\
&\overset{\text{Eq. (\ref{eq_sumf})}}{\scalebox{5}[1]{$\le$}} \frac{1}{\frac{1}{2}}\frac{1}{2 \varepsilon}d(x,y)+\frac{1}{\frac{1}{2}}\frac{1}{\frac{1}{2}}\sum_{w\in V}\lvert f_{w,2 \varepsilon}(x)-f_{w, 2 \varepsilon}(y)\rvert\\
&\overset{(*)}{\scalebox{2}[1]{$\le$}} \frac{1}{\varepsilon}d(x,y)+4N \frac{1}{2 \varepsilon}d(x,y)=\frac{2N+1}{\varepsilon}d(x,y),
\end{align*}
where in ($*$), we use the fact that $\lvert f_{w,2 \varepsilon}(x)-f_{w,2 \varepsilon}(y)\rvert\ne 0$ implies $d(x,w)<4 \varepsilon$. Hence, \ref{item_COLip} holds.
\end{proof}

The property of absolute continuity is preserved under linear combinations and under $\mathcal{E}$-convergence, as follows (see \cite[LEMMA 3.6 (a) and LEMMA 3.7 (a)]{KM20} for the Dirichlet form setting). The proof follows directly from the triangle inequality for $\Gamma(\cdot)(A)^{1/p}$ for any $A\in \mathcal{B}(X)$, and is therefore omitted.

\begin{lemma}\label{lem_AC_prep}
\hspace{0em}
\begin{enumerate}[label=(\arabic*),ref=(\arabic*)]
\item\label{lem_AC_prep1} If $f$, $g\in \mathcal{F}$ satisfy that $\Gamma(f)\ll m$ and $\Gamma(g)\ll m$, then for any $a,b\in \mathbb{R}$, we have $\Gamma(af+bg)\ll m$.
\item\label{lem_AC_prep2} If $\{f_n\}\subseteq \mathcal{F}$ and $f\in \mathcal{F}$ satisfy that $\Gamma(f_n)\ll m$ for any $n$, and $\lim_{n\to+\infty} \mathcal{E}(f_n-f)=0$, then $\Gamma(f)\ll m$.
\end{enumerate} 
\end{lemma}

\begin{proposition}[Energy dominance of $m$]\label{prop_dominant}
Assume \ref{eq_VD}, \ref{eq_PI}, \ref{eq_CS} and Equation (\ref{eq_Psi_AC}). Then $m$ is an energy-dominant measure of $(\mathcal{E},\mathcal{F})$, that is, $\Gamma(f)\ll m$ for any $f\in \mathcal{F}$.
\end{proposition}

\begin{proof}
Since $\mathcal{F}\cap C_c(X)$ is $\mathcal{E}_1$-dense in $\mathcal{F}$, by Lemma \ref{lem_AC_prep} \ref{lem_AC_prep2}, we only need to show that $\Gamma(f)\ll m$ for any $f\in \mathcal{F}\cap C_c(X)$.

By assumption, Lemma \ref{lem_Psi_p} holds, let $r_1\in(0,\mathrm{diam}(X))$ be the constant appearing in Equation (\ref{eq_Psi_p}). For any positive integer $n$ with $\frac{1}{n}<\frac{r_1}{2}$, let $V_n$ be a $\frac{1}{n}$-net, $\{\psi_z:z\in V_n\}\subseteq \mathcal{F}\cap C_c(X)$ the family of functions given by Lemma \ref{lem_partition}, and $f_n=\sum_{z\in V_n}f_{B(z,\frac{1}{n})}\psi_z$. Since $f\in C_c(X)$, we have $f_n$ is a finite linear combination of $\{\psi_z:z\in V_n\}$, which implies $f_n\in \mathcal{F}\cap C_c(X)$. By \ref{item_COAC} and Lemma \ref{lem_AC_prep} \ref{lem_AC_prep1}, we have $\Gamma(f_n)\ll m$.

We claim that $\{f_n\}$ converges uniformly to $f$, $\{f_n\}$ is $L^p$-convergent to $f$, and $\{f_n\}$ is $\mathcal{E}$-bounded. Indeed, for any $x\in X$, we have
\begin{align*}
&\lvert f_n(x)-f(x)\rvert\overset{\text{\ref{item_COunit}}}{\scalebox{4}[1]{$=$}}\lvert \sum_{z\in V_n}f_{B(z,\frac{1}{n})}\psi_z(x)-\sum_{z\in V_n}f(x)\psi_z(x)\rvert\le  \sum_{z\in V_n}\lvert f_{B(z,\frac{1}{n})}-f(x)\rvert\psi_z(x)\\
&\overset{\text{\ref{item_COspt}}}{\scalebox{4}[1]{$=$}}\sum_{z\in V_n:d(x,z)<\frac{2}{n}}\lvert f_{B(z,\frac{1}{n})}-f(x)\rvert\psi_z(x)\\
&\overset{f\in C_c(X)}{\scalebox{6}[1]{$\le$}}\sum_{z\in V_n:d(x,z)<\frac{2}{n}}\left(\sup\{|f(x_1)-f(x_2)|:d(x_1,x_2)<\frac{3}{n}\}\right)\psi_z(x)\\
&\overset{\text{\ref{item_COunit}}}{\scalebox{4}[1]{$\le$}}\sup\{|f(x_1)-f(x_2)|:d(x_1,x_2)<\frac{3}{n}\},
\end{align*}
hence
$$\sup_{x\in X}\lvert f_n(x)-f(x)\rvert\le\sup\{|f(x_1)-f(x_2)|:d(x_1,x_2)<\frac{3}{n}\}\to 0$$
as $n\to+\infty$, where we use the fact that $f\in C_c(X)$ is uniformly continuous. Hence, $\{f_n\}$ converges uniformly to $f$. Moreover, let $B(x_0,R)$ be a ball containing $\mathrm{supp}(f)$, then $\mathrm{supp}(f_n)\subseteq B(x_0,R+r_1)$ for any $n$, hence
\begin{align*}
\int_X|f_n-f|^p \mathrm{d}m\le \left(\sup_{x\in X}\lvert f_n(x)-f(x)\rvert\right)^pV(x_0,R+r_1)\to0
\end{align*}
as $n\to+\infty$, which gives $\{f_n\}$ is $L^p$-convergent to $f$.

For any $n$, for any $w\in V_n$, we have
\begin{align*}
&\Gamma(f_n)(B(w,\frac{1}{n}))\overset{\text{\ref{item_COunit}}}{\scalebox{4}[1]{$=$}}\Gamma\left(\sum_{z\in V_n}(f_{B(z,\frac{1}{n})}-f_{B(w,\frac{1}{n})})\psi_z+f_{B(w,\frac{1}{n})}\right)(B(w,\frac{1}{n}))\\
&\overset{\Gamma:\text{ strongly local}}{\underset{\text{\ref{item_COspt}}}{\scalebox{9}[1]{$=$}}}\Gamma\left(\sum_{z\in V_n:d(z,w)<\frac{3}{n}}(f_{B(z,\frac{1}{n})}-f_{B(w,\frac{1}{n})})\psi_z\right)(B(w,\frac{1}{n}))\\
&\le\left(\# \left\{z\in V_n:d(z,w)<\frac{3}{n}\right\}\right)^{p-1}\sum_{z\in V_n:d(z,w)<\frac{3}{n}}\lvert f_{B(z,\frac{1}{n})}-f_{B(w,\frac{1}{n})}\rvert^p \mathcal{E}(\psi_z),
\end{align*}
where we use the triangle inequality and H\"older's inequality in the last inequality. By \ref{eq_VD}, there exists some positive integer $N$ depending only on $C_{VD}$ such that $\# \{z\in V_n:d(z,w)<\frac{3}{n}\}\le N$. By \ref{item_COenergy}, we have
$$\mathcal{E}(\psi_z)\le C_1 \frac{V(z,\frac{1}{n})}{\left(\frac{1}{n}\right)^p},$$
where $C_1$ is the constant appearing therein. By \ref{eq_PI}, we have
\begin{align*}
&\lvert f_{B(z,\frac{1}{n})}-f_{B(w,\frac{1}{n})}\rvert^p\le\dashint_{B(z,\frac{1}{n})}\dashint_{B(w,\frac{1}{n})}|f(x)-f(y)|^pm(\mathrm{d}x)m(\mathrm{d}y)\\
&\overset{d(z,w)<\frac{3}{n}}{\scalebox{8}[1]{$\le$}}\frac{1}{V(z,\frac{1}{n})V(w,\frac{1}{n})}\int_{B(w,\frac{4}{n})}\int_{B(w,\frac{4}{n})}|f(x)-f(y)|^pm(\mathrm{d}x)m(\mathrm{d}y)\\
&\le\frac{2^pV(w,\frac{4}{n})}{V(z,\frac{1}{n})V(w,\frac{1}{n})}\int_{B(w,\frac{4}{n})}|f-f_{B(w,\frac{4}{n})}|^p \mathrm{d}m\le\frac{2^pC_{PI}V(w,\frac{4}{n})}{V(z,\frac{1}{n})V(w,\frac{1}{n})}\Psi(\frac{4}{n})\Gamma(f)(B(w,\frac{4A_{PI}}{n}))\\
&\overset{\text{\ref{eq_VD}, Eq. (\ref{eq_Psi_p})}}{\underset{\Psi:\text{ doubling}}{\scalebox{10}[1]{$\lesssim$}}}\frac{1}{n^pV(z,\frac{1}{n})}\Gamma(f)(B(w,\frac{4A_{PI}}{n})).
\end{align*}
Hence
\begin{align*}
&\Gamma(f_n)(B(w,\frac{1}{n}))\\
&\lesssim\sum_{z\in V_n:d(z,w)<\frac{3}{n}}\frac{1}{n^pV(z,\frac{1}{n})}\Gamma(f)(B(w,\frac{4A_{PI}}{n}))\frac{V(z,\frac{1}{n})}{\left(\frac{1}{n}\right)^p}\\
&\lesssim\Gamma(f)(B(w,\frac{4A_{PI}}{n})),
\end{align*}
which gives
$$\mathcal{E}(f_n)\le\sum_{w\in V_n}\Gamma(f_n)(B(w,\frac{1}{n}))\lesssim\sum_{w\in V_n}\Gamma(f)(B(w,\frac{4A_{PI}}{n}))=\int_{X}\left(\sum_{w\in V_n}1_{B(w,\frac{4A_{PI}}{n})}\right)\mathrm{d}\Gamma(f).$$
By \ref{eq_VD}, there exists some positive integer $M$ depending only on $C_{VD}$, $A_{PI}$ such that
$$\sum_{w\in V_n}1_{B(w,\frac{4A_{PI}}{n})}\le M1_{\cup_{w\in V_n}B(w,\frac{4A_{PI}}{n})},$$
hence $\mathcal{E}(f_n)\lesssim \mathcal{E}(f)$ for any $n$, $\{f_n\}$ is $\mathcal{E}$-bounded, which gives $\{f_n\}$ is $\mathcal{E}_1$-bounded.

Since $(\mathcal{F},\mathcal{E}_1)$ is a reflexive Banach space, by the Banach–Alaoglu theorem (see \cite[Theorem 3 in Chapter 12]{Lax02}), there exists a subsequence, still denoted by $\{f_n\}$, which is $\mathcal{E}_1$-weakly-convergent to some element $g\in \mathcal{F}$. By Mazur's lemma, here we refer to the version in \cite[Theorem 2 in Section V.1]{Yos95}, for any $n$, there exist $I_n\ge n$, $\lambda_k^{(n)}\ge0$ for $k=n,\ldots,I_n$ with $\sum_{k=n}^{I_n}\lambda_k^{(n)}=1$ such that $\{\sum_{k=n}^{I_n}\lambda_k^{(n)}f_k\}_{n}$ is $\mathcal{E}_1$-convergent to $g$, hence also $L^p$-convergent to $g$. Since $\{f_n\}$ is $L^p$-convergent to $f$, we have
$$\lVert {\sum_{k=n}^{I_n}\lambda_k^{(n)}f_k-f}\rVert_{L^p(X;m)}\le\sum_{k=n}^{I_n}\lambda_k^{(n)}\lVert {f_k-f}\rVert_{L^p(X;m)}\le\sup_{k\ge n}\lVert {f_k-f}\rVert_{L^p(X;m)}\to0$$
as $n\to+\infty$, which gives $f=g$. Hence $\{\sum_{k=n}^{I_n}\lambda_k^{(n)}f_k\}_{n}$ is $\mathcal{E}_1$-convergent to $f$. By Lemma \ref{lem_AC_prep} \ref{lem_AC_prep1}, we have $\Gamma(\sum_{k=n}^{I_n}\lambda_k^{(n)}f_k)\ll m$ for any $n$. By Lemma \ref{lem_AC_prep} \ref{lem_AC_prep2}, we have $\Gamma(f)\ll m$.
\end{proof}

\begin{proposition}[Minimality of $m$]\label{prop_minimal}
Assume \ref{eq_VD}, \ref{eq_PI}, \ref{eq_CS} and Equation (\ref{eq_Psi_AC}). If $\nu$ is an energy-dominant measure of $(\mathcal{E},\mathcal{F})$, that is, $\Gamma(f)\ll\nu$ for any $f\in \mathcal{F}$, then $m\ll\nu$.
\end{proposition}

\begin{proof}
Let $m=m_a+m_s$ be the Lebesgue decomposition of $m$ with respect to $\nu$, where $m_a\ll\nu$ and $m_s\perp\nu$. We only need to show that $m_s(X)=0$. We claim that there exist $C\in(0,+\infty)$, $R\in(0,\mathrm{diam}(X))$ such that for any $x\in X$, for any $r\in(0,R)$, we have
\begin{equation}\label{eq_min_mma}
m(B(x,r))\le Cm_a(B(x,r)).
\end{equation}
Then suppose $m_s(X)>0$, by the regularity of $m_s$, there exists a compact subset $K\subseteq X$ such that $m_s(K)>0$ and $m_a(K)=0$. For any $\varepsilon\in(0,R)$, let $V_{2 \varepsilon}$ be a $(2 \varepsilon)$-net of $(K,d)$. Since $K$ is compact, we have $V_{2 \varepsilon}$ is a finite set, which follows that
\begin{align*}
&0<m_s(K)=m(K)\le\sum_{z\in V_{2 \varepsilon}}m(B(z,2 \varepsilon))\\
&\overset{\text{\ref{eq_VD}}}{\scalebox{4}[1]{$\le$}}C_{VD}\sum_{z\in V_{2 \varepsilon}}m(B(z,\varepsilon))\overset{\text{Eq. (\ref{eq_min_mma})}}{\scalebox{6}[1]{$\le$}}C_{VD}C\sum_{z\in V_{2 \varepsilon}}m_a(B(z,\varepsilon))\\
&\overset{V_{2 \varepsilon}\text{: $(2 \varepsilon)$-net}}{\scalebox{8}[1]{$=$}}C_{VD}Cm_a \left(\bigcup_{z\in V_{2 \varepsilon}}B(z,\varepsilon)\right)\overset{V_{2 \varepsilon}\subseteq K}{\scalebox{5}[1]{$\le$}}C_{VD}Cm_a(K_\varepsilon),
\end{align*}
where $K_\varepsilon=\cup_{z\in K}B(z,\varepsilon)$. Since $K$ is compact, we have $\cap_{\varepsilon\in(0,R)}K_\varepsilon=K$. By the regularity of $m_a$, we have
$$0<m_s(K)\le C_{VD}C\lim_{\varepsilon\downarrow0}m_a(K_\varepsilon)=C_{VD}Cm_a(K)=0,$$
which gives a contradiction. Therefore, $m_s(X)=0$, $m=m_a\ll\nu$.

We only need to prove Equation (\ref{eq_min_mma}). By assumption, Lemma \ref{lem_Psi_p} holds, let $r_1\in(0,\mathrm{diam}(X))$, $C_1$ be the constants appearing therein, and Lemma \ref{lem_fxr} holds, let $C_2$ be the constant appearing therein. For any $x\in X$, for any $r\in(0,r_1)$, let $f_{x,r}=(1-\frac{d(x,\cdot)}{r})_+\in \mathcal{F}$ be the function given by Lemma \ref{lem_fxr}, then $\Gamma(f_{x,r})\le \frac{C_2^p}{r^p}m$. Since $m=m_a+m_s$, $m_a\ll \nu$, $m_s\perp\nu$, there exist disjoint measurable sets $E_1$, $E_2$ with $X=E_1\cup E_2$ such that $m_s(E_1)=0$ and $m_a(E_2)=\nu(E_2)=0$. Since $\Gamma(f_{x,r})\ll\nu$, we have $\Gamma(f_{x,r})(E_2)=0$. Then for any measurable set $U$, we have
$$\Gamma(f_{x,r})(U)=\Gamma(f_{x,r})(U\cap E_1)\le \frac{C_2^p}{r^p}m(U\cap E_1)=\frac{C_2^p}{r^p}m_a(U\cap E_1)=\frac{C_2^p}{r^p}m_a(U),$$
hence
\begin{equation}\label{eq_Gamma_ma}
\Gamma(f_{x,r})\le \frac{C_2^p}{r^p}m_a.
\end{equation}

By \ref{eq_PI}, we have
$$\int_{B(x,r)}\lvert f_{x,r}-(f_{x,r})_{B(x,r)}\rvert^p \mathrm{d}m\le C_{PI}\Psi(r)\int_{B(x,A_{PI}r)}\mathrm{d}\Gamma(f_{x,r}).$$
Since $f_{x,r}(y)\in[0,1]$ for any $y\in X$, we have $(f_{x,r})_{B(x,r)}\in[0,1]$. If $(f_{x,r})_{B(x,r)}\in[0,\frac{1}{2}]$, then since $f_{x,r}\ge \frac{3}{4}$ in $B(x,\frac{r}{4})$, we have
\begin{align*}
\int_{B(x,r)}|f_{x,r}-(f_{x,r})_{B(x,r)}|^p \mathrm{d}m\ge \frac{1}{4^p}m(B(x,\frac{r}{4}))\overset{\text{\ref{eq_VD}}}{\scalebox{4}[1]{$\ge$}}\frac{1}{4^pC_{VD}^2}m(B(x,r)).
\end{align*}
If $(f_{x,r})_{B(x,r)}\in[\frac{1}{2},1]$, then since $f_{x,r}\le \frac{1}{4}$ in $B(x,r)\backslash B(x,\frac{3r}{4})$, we have
\begin{align*}
\int_{B(x,r)}|f_{x,r}-(f_{x,r})_{B(x,r)}|^p \mathrm{d}m\ge \frac{1}{4^p}m(B(x,r)\backslash B(x,\frac{3r}{4})).
\end{align*}
By \ref{eq_CC}, there exists a ball $B(y,\frac{r}{16})\subseteq B(x,r)\backslash B(x,\frac{3r}{4})$, hence
$$\int_{B(x,r)}|f_{x,r}-(f_{x,r})_{B(x,r)}|^p \mathrm{d}m\ge \frac{1}{4^p}m(B(y,\frac{r}{16}))\overset{\text{\ref{eq_VD}}}{\scalebox{4}[1]{$\ge$}}\frac{1}{4^pC_{VD}^6}m(B(x,r)).$$
Therefore
\begin{align*}
&m(B(x,r))\le4^pC_{VD}^6\int_{B(x,r)}|f_{x,r}-(f_{x,r})_{B(x,r)}|^p \mathrm{d}m\\
&\le4^pC_{VD}^6C_{PI}\Psi(r)\int_{B(x,A_{PI}r)}\mathrm{d}\Gamma(f_{x,r})\\
&\overset{\text{Eq. (\ref{eq_Gamma_ma})}}{\underset{\text{Eq. (\ref{eq_Psi_p})}}{\scalebox{6}[1]{$\le$}}}4^pC_{VD}^6C_{PI}C_1r^p\frac{C_2^p}{r^p}m_a({B(x,A_{PI}r)})\\
&=4^pC_1C_2^pC_{PI}C_{VD}^6m_a({B(x,A_{PI}r)}),
\end{align*}
which gives
$$m(B(x,A_{PI}r))\overset{\text{\ref{eq_VD}}}{\scalebox{4}[1]{$\le$}}C_{VD}^{\log_2A_{PI}+1}m(B(x,r))\le Cm_a({B(x,A_{PI}r)}),$$
where $C=4^pC_1C_2^pC_{PI}C_{VD}^{\log_2A_{PI}+7}$. Therefore, we have Equation (\ref{eq_min_mma}) holds with $R=\min\{A_{PI}r_1,\mathrm{diam}(X)\}$.
\end{proof}

\begin{proposition}\label{prop_biLip}
Assume \ref{eq_VD}, \ref{eq_PI}, \ref{eq_CS} and Equation (\ref{eq_Psi_AC}). Then $\rho$ is a geodesic metric on $X$, and $\rho$ is bi-Lipschitz equivalent to $d$.
\end{proposition}

\begin{proof}
By assumption, Lemma \ref{lem_Psi_p} holds, let $r_1\in(0,\mathrm{diam}(X))$, $C_1$ be the constants appearing therein, and Lemma \ref{lem_fxr} holds, let $C_2$ be the constant appearing therein.

For any $x$, $y\in X$, for any $r\in(0,r_1)$, let $f_{x,r}=(1-\frac{d(x,\cdot)}{r})_+\in \mathcal{F}$ be given by Lemma \ref{lem_fxr}, then $\Gamma(f_{x,r})\le \frac{C_2^p}{r^p}m$, that is, $\Gamma(\frac{r}{C_2}f_{x,r})\le m$. If $d(x,y)<r_1$, then for any $r\in(d(x,y),r_1)$, we have $\rho(x,y)\ge \frac{r}{C_2}f_{x,r}(x)-\frac{r}{C_2}f_{x,r}(y)=\frac{1}{C_2}d(x,y)$; if $d(x,y)\ge r_1$, then for any $r\in(0,r_1)$, we have $\rho(x,y)\ge\frac{r}{C_2}f_{x,r}(x)-\frac{r}{C_2}f_{x,r}(y)= \frac{r}{C_2}$, letting $r\uparrow r_1$, we have $\rho(x,y)\ge \frac{r_1}{C_2}$, or equivalently, if $\rho(x,y)<\frac{r_1}{C_2}$, then $d(x,y)<r_1$.

On the other hand, for any $x$, $y\in X$, for any $f\in \mathcal{F}_{\mathrm{loc}}\cap C(X)$ with $\Gamma(f)\le m$, by \cite[Lemma 3.4]{Yan25a}, we have
$$\lvert f(x)-f(y)\rvert^p\le 2C_3\Psi(d(x,y)),$$
where $C_3$ is the constant appearing therein, hence $\rho(x,y)\le(2C_3)^{1/p}\Psi(d(x,y))^{1/p}<+\infty$. In particular, if $d(x,y)<r_1$, then by Lemma \ref{lem_Psi_p}, we have $\rho(x,y)\le (2C_1C_3)^{1/p}d(x,y)$.

In summary, we have
\begin{equation}\label{eq_biLip1}
\rho(x,y)<+\infty\text{ for any }x,y\in X,
\end{equation}
%\begin{equation}\label{eq_biLip2}
%d(x,y)<r_1\text{ for any }x,y\in X\text{ with }\rho(x,y)<\frac{r_1}{C_2},
%\end{equation}
\begin{equation}\label{eq_biLip2}
\frac{1}{C_4}d(x,y)\le \rho(x,y)\le C_4d(x,y)\text{ for any }x,y\in X\text{ with }d(x,y)<r_1\text{ or }\rho(x,y)<\frac{r_1}{C_2},
\end{equation}
with $C_4=\max\{C_2,(2C_1C_3)^{1/p}\}$. If $\rho(x,y)=0$, then by Equation (\ref{eq_biLip2}), we have $d(x,y)=0$, hence $x=y$. Combining this with Equation (\ref{eq_biLip1}), we have $\rho$ is a metric. By Equation (\ref{eq_biLip2}), \ref{asm_A} holds. Then by Proposition \ref{prop_geodesic}, we have $\rho$ is a geodesic metric.

For any $x$, $y\in X$. Firstly, take an integer $n\ge1$ such that $C_{cc}\frac{d(x,y)}{n}<r_1$, where $C_{cc}$ is the constant in \ref{eq_CC}, then there exists a sequence $\{x_k:0\le k\le n\}$ with $x_0=x$ and $x_n=y$ such that $d(x_k,x_{k-1})\le C_{cc}\frac{d(x,y)}{n}<r_1$ for any $k=1,\ldots,n$. By Equation (\ref{eq_biLip2}), we have $\rho(x_k,x_{k-1})\le C_4d(x_k,x_{k-1})$. Hence
$$\rho(x,y)\le\sum_{k=1}^n\rho(x_k,x_{k-1})\le C_4\sum_{k=1}^nd(x_k,x_{k-1})\le C_4C_{cc}d(x,y).$$
Secondly, take an integer $n\ge1$ such that $\frac{\rho(x,y)}{n}<\frac{r_1}{C_2}$. Since $\rho$ is a geodesic metric, there exists a sequence $\{y_k:0\le k\le n\}$ with $y_0=x$ and $y_n=y$ such that $\rho(y_k,y_{k-1})=\frac{\rho(x,y)}{n}<\frac{r_1}{C_2}$ for any $k=1,\ldots,n$. By Equation (\ref{eq_biLip2}), we have $d(y_k,y_{k-1})\le C_4\rho(y_k,y_{k-1})$. Hence
$$d(x,y)\le\sum_{k=1}^nd(y_k,y_{k-1})\le C_4\sum_{k=1}^n\rho(y_k,y_{k-1})=C_4\rho(x,y).$$
Therefore, $\rho$ is bi-Lipschitz equivalent to $d$.
\end{proof}

\begin{proof}[Proof of Theorem \ref{thm_AC}]
It follows directly from Proposition \ref{prop_dominant}, Proposition \ref{prop_minimal}, and Proposition \ref{prop_biLip}.
\end{proof}

\section{Proof of Theorem \ref{thm_dichotomy}}\label{sec_dichotomy}

For any $\alpha\in(0,+\infty)$, we have the following definition of Besov spaces:
\begin{align*}
&B^{p,\alpha}(X)\\
&=\left\{f\in L^p(X;m):\sup_{r\in(0,\mathrm{diam}(X))}\frac{1}{r^{p\alpha}}\int_X\dashint_{B(x,r)}|f(x)-f(y)|^pm(\mathrm{d} y)m(\mathrm{d} x)<+\infty\right\}.
\end{align*}
Obviously, $B^{p,\alpha}(X)$ is decreasing in $\alpha$ and may become trivial if $\alpha$ is too large. We define the following critical exponent
\begin{align*}
&\alpha_p(X)=\sup\left\{\alpha\in(0,+\infty):B^{p,\alpha}(X)\text{ contains non-constant functions}\right\}\le+\infty.
\end{align*}
Notably the value of $\alpha_p(X)$ depends \emph{only} on the metric measure space $(X,d,m)$. We have some basic properties of $\alpha_p(X)$ as follows.

\begin{lemma}[{\cite[Theorem 4.1]{Bau24}}]\label{lem_alpha_basic}
\hfill
\begin{enumerate}[label=(\roman*)]
\item For any $p\in(1,+\infty)$, we have $\alpha_p(X)\ge1$.
\item The function $p\mapsto p\alpha_p(X)$ is monotone increasing for $p\in(1,+\infty)$.
\item The function $p\mapsto \alpha_p(X)$ is monotone decreasing for $p\in(1,+\infty)$.
\end{enumerate}
Hence
\begin{enumerate}[label=(\alph*)]
\item For $p\in(1,+\infty)$, the functions $p\mapsto p\alpha_p(X)$ and $p\mapsto \alpha_p(X)$ are continuous.
\item If $\alpha_p(X)<+\infty$ for some $p\in(1,+\infty)$, then $\alpha_p(X)<+\infty$ for all $p\in(1,+\infty)$.
\end{enumerate}
\end{lemma}

The value of $\alpha_p(X)$ can be determined once certain functional inequalities are satisfied as follows.

\begin{lemma}[{\cite[Theorem 4.6]{Shi24a}}]\label{lem_alpha_beta}
Assume \ref{eq_VD}, \hyperlink{eq_PIbeta}{PI$_p$($\beta_p$)}, \hyperlink{eq_ucapbeta}{$\text{cap}_p(\beta_p)_{\le}$}. Then $\alpha_p(X)=\frac{\beta_p}{p}$.
\end{lemma}

For $x\in X$, for a function $u$ defined in an open neighborhood of $x$, its pointwise Lipschitz constant at $x$ is defined as
$$\mathrm{Lip}\ u(x)=\lim_{r\downarrow0}\sup_{y:d(x,y)\in(0,r)}\frac{\lvert u(x)-u(y)\rvert}{d(x,y)}.$$
We say a function $u$ defined in $X$ is Lipschitz if there exists $K\in(0,+\infty)$ such that $\lvert u(x)-u(y)\rvert\le Kd(x,y)$ for any $x$, $y\in X$. Let $\mathrm{Lip}(X)$ be the family of all Lipschitz functions.

\begin{proposition}\label{prop_PI_open}
Assume \ref{eq_VD}. Let $p\in(1,+\infty)$. If there exists a $p$-energy $(\mathcal{E},\mathcal{F})$ such that \hyperlink{eq_PIbeta}{PI$_p$($p$)}, \hyperlink{eq_CSbeta}{$\text{CS}_p(p)$} hold, then there exists $\varepsilon>0$ such that $\alpha_q(X)=1$ for any $q\in(p-\varepsilon,+\infty)$.
\end{proposition}

\begin{remark}
If we replace \hyperlink{eq_CSbeta}{$\text{CS}_p(p)$} by \hyperlink{eq_ucapbeta}{$\text{cap}_p(p)_{\le}$}, then by Lemma \ref{lem_alpha_beta}, we obtain $\alpha_p(X)=1$. Combining this with the monotonicity of $p\mapsto \alpha_p(X)\ge1$ from Lemma \ref{lem_alpha_basic}, it follows that $\alpha_q(X)=1$ for any $q\in[p,+\infty)$. The key point of our result is that if the stronger condition \hyperlink{eq_CSbeta}{$\text{CS}_p(p)$} holds, then this equality can be ``self-improved" to hold in a slightly larger open interval $(p-\varepsilon,+\infty)$.
\end{remark}

\begin{proof}[Proof of Proposition \ref{prop_PI_open}]
By Theorem \ref{thm_AC}, $m$ is a minimal energy-dominant measure of $(\mathcal{E},\mathcal{F})$, $\rho$ is a geodesic metric on $X$, and $\rho$ is bi-Lipschitz equivalent to $d$; let $C_1$ denote the Lipschitz constant associated with this equivalence. Notably, \ref{asm_A} holds. By \ref{eq_VD} for $d$, we have Equation (\ref{eq_rho_VD}) for $\rho$, then by Theorem \ref{thm_Lip}, we have $\mathrm{Lip}(X)=\mathrm{Lip}_\rho(X)\subseteq \mathcal{F}_{\mathrm{loc}}$, and for any $u\in \mathrm{Lip}(X)$, we have $\Gamma(u)\le (\mathrm{Lip}_\rho u)^p m\le C_1^p(\mathrm{Lip}\ u)^pm$. Hence for any ball $B(x_0,R)$, we have
\begin{align*}
&\dashint_{B(x_0,R)}\lvert u-u_{B(x_0,R)}\rvert \mathrm{d}m\le \left(\dashint_{B(x_0,R)}\lvert u-u_{B(x_0,R)}\rvert^p \mathrm{d}m\right)^{1/p}\\
&\overset{\text{\hyperlink{eq_PIbeta}{PI$_p$($p$)}}}{\scalebox{4}[1]{$\le$}} \left(\frac{1}{V(x_0,R)}C_{PI}R^p\int_{B(x_0,A_{PI}R)}\mathrm{d}\Gamma(u)\right)^{1/p} \\
&\overset{\text{\ref{eq_VD}}}{\scalebox{4}[1]{$\le$}} \left(\frac{C_2}{V(x_0,A_{PI}R)}C_{PI}R^p\int_{B(x_0,A_{PI}R)}C_1^p(\mathrm{Lip}\ u)^p\mathrm{d}m\right)^{1/p}\\
&=C_1(C_2C_{PI})^{1/p}R \left(\dashint_{B(x_0,A_{PI}R)}(\mathrm{Lip}\ u)^p\mathrm{d}m\right)^{1/p},
\end{align*}
where $C_2$ is some positive constant depending only on $C_{VD}$, $A_{PI}$. Let $C=C_1(C_2C_{PI})^{1/p}$, $A=A_{PI}$, then $(X,d,m)$ supports the following $(1,p)$-Poincar\'e inequality \hypertarget{eq_PI_1p}{$\mathrm{PI}_{\mathrm{Lip}}(1,p)$}: for any ball $B(x_0,R)$, for any $u\in \mathrm{Lip}(X)$, we have
\begin{equation*}
\dashint_{B(x_0,R)}|u-u_{B(x_0,R)}| \mathrm{d} m\le CR\left(\dashint_{B(x_0,AR)}(\mathrm{Lip}\ u)^p \mathrm{d}m\right)^{1/p}.
\end{equation*}
By \cite[THEOREM 1.0.1]{KZ08}, there exists $\varepsilon>0$ such that for any $q\in(p-\varepsilon,+\infty)$, $(X,d,m)$ supports a $(1,q)$-Poincar\'e inequality \hyperlink{eq_PI_1p}{$\mathrm{PI}_{\mathrm{Lip}}(1,q)$}. By \cite[Theorem 5.1, Remark 5.2]{Bau24}, the condition $\mathcal{P}(q,1)$ holds (see \cite[Definition 4.5]{Bau24} for its definition). Consequently, \cite[Lemma 4.7]{Bau24} yields $\alpha_q(X)=1$.
\end{proof}

\begin{proof}[Proof of Theorem \ref{thm_dichotomy}]
For any $p\in I$, by Lemma \ref{lem_alpha_beta}, we have $\alpha_p(X)=\frac{\beta_p}{p}<+\infty$. Let
$$J=\{p\in I:\alpha_p(X)=1\}.$$
We only need to show that either $J=\emptyset$ or $J=I$. Indeed, suppose that $J\ne\emptyset$ but $J\ne I$. By Lemma \ref{lem_alpha_basic}, we have $p\mapsto \alpha_p(X)$ is monotone decreasing and continuous, hence $J=[p,+\infty)\cap I$ is an interval for some $p\in I$. However, since $p\in J$, $\beta_p=p$, under \ref{eq_VD}, \hyperlink{eq_PIbeta}{PI$_p$($p$)}, \hyperlink{eq_CSbeta}{$\text{CS}_p(p)$}, by Proposition \ref{prop_PI_open}, there exists $\varepsilon>0$ such that $(p-\varepsilon,+\infty)\cap I\subseteq J$, contradiction.
\end{proof}

\bibliographystyle{plain}
%\bibliography{/Users/meng/Dropbox/myref}

\begin{thebibliography}{10}

\bibitem{AB25}
Patricia Alonso-Ruiz and Fabrice Baudoin.
\newblock Korevaar-{S}choen {$p$}-energies and their {$\Gamma $}-limits on
  {C}heeger spaces.
\newblock {\em Nonlinear Anal.}, 256:Paper No. 113779, 22, 2025.

\bibitem{AES25a}
Riku {Anttila}, Sylvester {Eriksson-Bique}, and Ryosuke {Shimizu}.
\newblock {Construction of self-similar energy forms and singularity of Sobolev
  spaces on Laakso-type fractal spaces}.
\newblock {\em arXiv e-prints}, page arXiv:2503.13258, March 2025.

\bibitem{BB89}
Martin~T. Barlow and Richard~F. Bass.
\newblock The construction of {B}rownian motion on the {S}ierpi\'{n}ski carpet.
\newblock {\em Ann. Inst. H. Poincar\'{e} Probab. Statist.}, 25(3):225--257,
  1989.

\bibitem{BB90}
Martin~T. Barlow and Richard~F. Bass.
\newblock On the resistance of the {S}ierpi\'{n}ski carpet.
\newblock {\em Proc. Roy. Soc. London Ser. A}, 431(1882):345--360, 1990.

\bibitem{BB92}
Martin~T. Barlow and Richard~F. Bass.
\newblock Transition densities for {B}rownian motion on the {S}ierpi\'{n}ski
  carpet.
\newblock {\em Probab. Theory Related Fields}, 91(3-4):307--330, 1992.

\bibitem{BBS90}
Martin~T. Barlow, Richard~F. Bass, and John~D. Sherwood.
\newblock Resistance and spectral dimension of {S}ierpi\'{n}ski carpets.
\newblock {\em J. Phys. A}, 23(6):L253--L258, 1990.

\bibitem{BP88}
Martin~T. Barlow and Edwin~A. Perkins.
\newblock Brownian motion on the {S}ierpi\'{n}ski gasket.
\newblock {\em Probab. Theory Related Fields}, 79(4):543--623, 1988.

\bibitem{Bau24}
Fabrice Baudoin.
\newblock Korevaar-{S}choen-{S}obolev spaces and critical exponents in metric
  measure spaces.
\newblock {\em Ann. Fenn. Math.}, 49(2):487--527, 2024.

\bibitem{BC23}
Fabrice Baudoin and Li~Chen.
\newblock Sobolev spaces and {P}oincar\'e{} inequalities on the {V}icsek
  fractal.
\newblock {\em Ann. Fenn. Math.}, 48(1):3--26, 2023.

\bibitem{BH99}
Martin~R. Bridson and Andr\'e Haefliger.
\newblock {\em Metric spaces of non-positive curvature}, volume 319 of {\em
  Grundlehren der mathematischen Wissenschaften [Fundamental Principles of
  Mathematical Sciences]}.
\newblock Springer-Verlag, Berlin, 1999.

\bibitem{CGQ22}
Shiping Cao, Qingsong Gu, and Hua Qiu.
\newblock {$p$}-energies on p.c.f. self-similar sets.
\newblock {\em Adv. Math.}, 405:Paper No. 108517, 58, 2022.

\bibitem{Che99}
Jeff Cheeger.
\newblock Differentiability of {L}ipschitz functions on metric measure spaces.
\newblock {\em Geom. Funct. Anal.}, 9(3):428--517, 1999.

\bibitem{CGYZ24}
Aobo {Chen}, Jin {Gao}, Zhenyu {Yu}, and Junda {Zhang}.
\newblock {Besov-Lipschitz norm and $p$-en-ergy measure on scale-irregular
  Vicsek sets}.
\newblock {\em J. Fractal Geom.}, 2025.
\newblock published online first.

\bibitem{FOT11}
Masatoshi Fukushima, Yoichi Oshima, and Masayoshi Takeda.
\newblock {\em Dirichlet forms and symmetric {M}arkov processes}, volume~19 of
  {\em De Gruyter Studies in Mathematics}.
\newblock Walter de Gruyter \& Co., Berlin, extended edition, 2011.

\bibitem{HKKZ00}
Ben~M. Hambly, Takashi Kumagai, Shigeo Kusuoka, and Xian~Yin Zhou.
\newblock Transition density estimates for diffusion processes on homogeneous
  random {S}ierpinski carpets.
\newblock {\em J. Math. Soc. Japan}, 52(2):373--408, 2000.

\bibitem{HPS04}
Paul~Edward Herman, Roberto Peirone, and Robert~S. Strichartz.
\newblock {$p$}-energy and {$p$}-harmonic functions on {S}ierpinski gasket type
  fractals.
\newblock {\em Potential Anal.}, 20(2):125--148, 2004.

\bibitem{Hin10}
Masanori Hino.
\newblock Energy measures and indices of {D}irichlet forms, with applications
  to derivatives on some fractals.
\newblock {\em Proc. Lond. Math. Soc. (3)}, 100(1):269--302, 2010.

\bibitem{KM20}
Naotaka Kajino and Mathav Murugan.
\newblock On singularity of energy measures for symmetric diffusions with full
  off-diagonal heat kernel estimates.
\newblock {\em Ann. Probab.}, 48(6):2920--2951, 2020.

\bibitem{KS24a}
Naotaka {Kajino} and Ryosuke {Shimizu}.
\newblock {Contraction properties and differentiability of $p$-energy forms
  with applications to nonlinear potential theory on self-similar sets}.
\newblock {\em arXiv e-prints}, page arXiv:2404.13668v2, April 2024.

\bibitem{KZ08}
Stephen Keith and Xiao Zhong.
\newblock The {P}oincar\'e{} inequality is an open ended condition.
\newblock {\em Ann. of Math. (2)}, 167(2):575--599, 2008.

\bibitem{Kig89}
Jun Kigami.
\newblock A harmonic calculus on the {S}ierpi\'{n}ski spaces.
\newblock {\em Japan J. Appl. Math.}, 6(2):259--290, 1989.

\bibitem{Kig23}
Jun Kigami.
\newblock {\em Conductive homogeneity of compact metric spaces and construction
  of {$p$}-energy}, volume~5 of {\em Memoirs of the European Mathematical
  Society}.
\newblock European Mathematical Society (EMS), Berlin, [2023] \copyright 2023.

\bibitem{KSZ12}
Pekka Koskela, Nageswari Shanmugalingam, and Yuan Zhou.
\newblock {$L^\infty$}-variational problem associated to {D}irichlet forms.
\newblock {\em Math. Res. Lett.}, 19(6):1263--1275, 2012.

\bibitem{KZ12}
Pekka Koskela and Yuan Zhou.
\newblock Geometry and analysis of {D}irichlet forms.
\newblock {\em Adv. Math.}, 231(5):2755--2801, 2012.

\bibitem{KZ92}
Shigeo Kusuoka and Xianyin Zhou.
\newblock Dirichlet forms on fractals: {P}oincar\'{e} constant and resistance.
\newblock {\em Probab. Theory Related Fields}, 93(2):169--196, 1992.

\bibitem{Lax02}
Peter~D. Lax.
\newblock {\em Functional analysis}.
\newblock Pure and Applied Mathematics (New York). Wiley-Interscience [John
  Wiley \& Sons], New York, 2002.

\bibitem{MS23}
Mathav Murugan and Ryosuke Shimizu.
\newblock First-order {S}obolev spaces, self-similar energies and energy
  measures on the {S}ierpi\'nski carpet.
\newblock {\em Comm. Pure Appl. Math.}, 78(9):1523--1608, 2025.

\bibitem{Nak85}
Shintaro Nakao.
\newblock Stochastic calculus for continuous additive functionals of zero
  energy.
\newblock {\em Z. Wahrsch. Verw. Gebiete}, 68(4):557--578, 1985.

\bibitem{Sas25}
K{\^o}hei {Sasaya}.
\newblock {Construction of $p$-energy measures associated with strongly local
  $p$-energy forms}.
\newblock {\em arXiv e-prints}, page arXiv:2502.10369v3, February 2025.

\bibitem{Shi24}
Ryosuke Shimizu.
\newblock Construction of {$p$}-energy and associated energy measures on
  {S}ierpi\'{n}ski carpets.
\newblock {\em Trans. Amer. Math. Soc.}, 377(2):951--1032, 2024.

\bibitem{Shi24a}
Ryosuke {Shimizu}.
\newblock Characterizations of {S}obolev {F}unctions {V}ia {B}esov-type
  {E}nergy {F}unctionals in {F}ractals.
\newblock {\em Potential Anal.}, 63(4):2121–2156, 2025.

\bibitem{Stu94}
Karl-Theodor Sturm.
\newblock Analysis on local {D}irichlet spaces. {I}. {R}ecurrence,
  conservativeness and {$L^p$}-{L}iouville properties.
\newblock {\em J. Reine Angew. Math.}, 456:173--196, 1994.

\bibitem{Stu95}
Karl-Theodor Sturm.
\newblock On the geometry defined by {D}irichlet forms.
\newblock In {\em Seminar on {S}tochastic {A}nalysis, {R}andom {F}ields and
  {A}pplications ({A}scona, 1993)}, volume~36 of {\em Progr. Probab.}, pages
  231--242. Birkh\"auser, Basel, 1995.

\bibitem{Yan25d}
Meng {Yang}.
\newblock {Energy inequalities for cutoff functions of $p$-energies on metric
  measures spaces}.
\newblock {\em arXiv e-prints}, page arXiv:2507.08577v3, July 2025.

\bibitem{Yan25c}
Meng {Yang}.
\newblock {On singularity of $p$-energy measures on metric measure spaces}.
\newblock {\em arXiv e-prints}, page arXiv:2505.12468v1, May 2025.

\bibitem{Yan25a}
Meng {Yang}.
\newblock {$p$-Poincar\'e inequalities and cutoff Sobolev inequalities on
  metric measure spaces}.
\newblock {\em arXiv e-prints}, page arXiv:2504.09503v2, April 2025.

\bibitem{Yos95}
K{\=o}saku Yosida.
\newblock {\em Functional analysis}.
\newblock Classics in Mathematics. Springer-Verlag, Berlin, 1995.
\newblock Reprint of the sixth (1980) edition.

\end{thebibliography}

\end{document}